\documentclass[twocolumn]{autart}

\usepackage{graphics}
\usepackage{graphicx}
\usepackage{amssymb}
\usepackage{verbatim}
\usepackage{amsxtra}
\usepackage{epsfig}
\usepackage{subfig}
\usepackage{float}
\usepackage[breaklinks = true]{hyperref} 
\usepackage{amsmath}
\usepackage{accents}

\usepackage{mathtools}

\usepackage{latexsym}
\usepackage{ifthen}
\usepackage{color}
\usepackage{hyperref}
\usepackage{balance}
\usepackage{subeqnarray}
\usepackage{caption}
\usepackage[dvipsnames]{xcolor}

\usepackage{rgsEnvironments_auto}
\usepackage{rgsMacros_auto}

\usepackage{enumitem}
\usepackage{algorithm,algpseudocode}

\usepackage{siunitx}

\usepackage{optidef}

\usepackage{pgfplots}

\usetikzlibrary{external}
\usepgfplotslibrary{external}
\usepackage{pgfplots}
\pgfplotsset{compat=newest,every axis/.append style={
                    label style={font=\footnotesize},
                    tick label style={font=\footnotesize}  
                    }}

\usepgfplotslibrary{fillbetween}

\DeclareMathOperator*{\interior}{int}
\DeclareMathOperator*{\closure}{cl}

\usepackage{amsmath}

\DeclareMathOperator*{\argmin}{arg\,min}

\newcommand{\Anorm}[1]{\left|#1\right|_{\mathcal{A}}}

\newcommand{\T}{\mathcal{T}}

\graphicspath{{./Figures/}}

\edef\endfrontmatter{%
  \unexpanded\expandafter{\endfrontmatter}
  \noexpand\endNoHyper 
}

\begin{document}


\begin{frontmatter}
\runtitle{Existence and Regularity of Optimal Controls for Hybrid Dynamical Systems}  

\title{Regularity of Optimal Solutions and the Optimal Cost for Hybrid Dynamical Systems via Reachability Analysis\thanksref{footnoteinfo}} 

\thanks[footnoteinfo]{This paper was not presented at any IFAC 
meeting. Corresponding author B.~Alt{\i}n. Research partially supported by the National Science Foundation under Grant no. ECS-1710621, Grant no. CNS-1544396, and Grant no. CNS-2039054, by the Air Force Office of Scientific Research under Grant no. FA9550-19-1-0053, Grant no. FA9550-19-1-0169, and Grant no. FA9550-20-1-0238, and by the Army Research Office under Grant no. W911NF-20-1-0253.}

\author[berk]{Berk Alt{\i}n}\ead{berkaltin@ucsc.edu},    
\author[berk]{Ricardo G. Sanfelice}\ead{ricardo@ucsc.edu}  

\address[berk]{Department of Electrical and Computer Engineering, University of California, Santa Cruz}


\begin{abstract}
For a general optimal control problem for dynamical systems with hybrid dynamics, we study the dependency  of the optimal cost and of the value function on the initial conditions, parameters, and perturbations. We show that upper and lower semicontinuous dependence of solutions on initial conditions -- properties that are captured by outer and inner well-posedness, respectively -- lead to the existence of a solution to the hybrid optimal control problem and upper/lower semicontinuity of the optimal cost. In particular, by exploiting properties of finite horizon reachable sets for hybrid systems, we show that the optimal cost varies upper semicontinuously when the hybrid system is (nominally) outer well-posed, and lower semicontinuously when it is (nominally) inner well-posed and an additional assumption requiring partial knowledge of solutions. Consequently, when the system is both (nominally) inner and outer well-posed and the aforementioned assumption holds, the optimal cost varies continuously, and optimal solutions vary upper semicontinuously. We further show that even in the absence of this solution-based assumption, the optimal cost can be continuously approximated.  The results are demonstrated in finite horizon optimization problems with hybrid dynamics, both theoretically and numerically.

\end{abstract}

\end{frontmatter}


\section{Introduction}
\label{sec:intro}
Models and algorithms characterized by the interplay of 
continuous-time dynamics and instantaneous changes 
have become prevalent due to their capabilities of 
leading to solutions to control problems that classical techniques
cannot solve, or simply do not apply.
Examples of outstanding control problems
that such {\em hybrid techniques} have been able to 
solve include 
the design of event-triggered control algorithms
\cite{Tabuada.07,Postoyan.ea.TAC.15,230},
stabilization over networks \cite{NesicTeel04a,Hespanha.ea.07}, 
observers and synchronization strategies under intermittent 
information \cite{136,190},
and control of mechanical systems exhibiting impacts \cite{RonsseLefevreSepulchre07TR,79}.
These advances have been enabled by 
the modeling, analysis, and design techniques 
for {\em hybrid dynamical systems}.
A hybrid dynamical system, or just a {\em hybrid system}, is a dynamical system that exhibits characteristics of both continuous-time and discrete-time dynamical systems.

Numerous tools are available in the literature for the study of hybrid systems, 
in particular, for hybrid systems modeled as
hybrid automata \cite{LygerosJohanssonSimicZhangSastry03,BranickyBorkarMitter98,vanderSchaftSchumacher00},
impulsive systems \cite{HaddadChellaboinaNersesov06,AubinLygerosQuincampoixSastrySeube02},
and 
hybrid inclusions \cite{hybridbook,220}.
The literature is rich in tools for the analysis of 
reachability \cite{LygerosTomlinSastry99,cdc2020,semantics},
asymptotic stability \cite{LygerosJohanssonSimicZhangSastry03,HaddadChellaboinaNersesov06,hybridbook},
forward invariance \cite{AubinLygerosQuincampoixSastrySeube02,185,9079562},
control design \cite{220},
and robustness \cite{hybridbook,220}.
On the other hand, 
optimality for hybrid systems is much less mature.
Initial results on optimality of trajectories
 over finite horizons were developed in \cite{Sussmann99}, 
 including a maximum principle for optimality, for a class of switched systems.
This result was extended in 
\cite{ShaikhCaines07,Pakniyat.Caines.20.TAC} to a broader class of systems, one allowing for state resets -- 
the models considered are in the spirit of hybrid automata.
More recently, 
linear-quadratic control for a class of hybrid systems with a sample-and-hold structure was considered in \cite{Possieri.Teel.16.CDC,Cristofaro.ea.18.ECC}.
In particular, the development in \cite{Possieri.Teel.16.CDC} is within the hybrid inclusions framework
in  \cite{hybridbook,220}, for the special case when
the continuous dynamics are modeled by a differential equation that is linear and the discrete dynamics are governed by a linear difference equation.
The problem of guaranteeing existence of 
optimal control inputs 
for a class of hybrid systems 
was studied in \cite{goebel2019}. 
The hybrid inclusions framework is employed in \cite{goebel2019} 
and the conditions for existence of optimal control inputs
require the continuous dynamics of the
system to be governed by a differential equation whose right-hand side
is affine in the control input.
Optimality of static state-feedback laws for 
hybrid inclusions with continuous and discrete dynamics modeled by 
(single-valued) nonlinear maps was studied in \cite{211}.
Infinitesimal conditions involving a Lyapunov-like function
are presented in \cite{211} to guarantee optimality over the infinite (hybrid) horizon.
The finite horizon optimization problem for the same broad class of hybrid systems
was formulated and developed in a sequence of papers leading to a model predictive control
framework; see \cite{ALTIN2018128,acc,9029729,9147972}.

Though the advances cited above have contributed to optimal control for hybrid systems, some of the key properties of the optimal control problem associated to general hybrid systems, wherein trajectories are constrained to evolve continuously ({\em flow}) in certain regions of the state space and to exhibit instantaneous changes ({\em jump}) under certain conditions, have not been yet revealed in the literature. Specifically, the regularity properties of the optimal cost, in particular, (semi) continuous dependence of the optimal cost and optimal trajectories on the constraints on where the trajectories can flow or jump have not yet been investigated. Very importantly, conditions enabling the approximations of the optimal cost in a continuous manner are not available in the literature. Indeed, results that permit relating the effect of varying parameters and initial conditions when they approach nominal values, the expectation being that the optimal cost also approaches its nominal value, are missing. Understanding such  a dependency is critical due to the fact that it is unavoidable to numerically compute trajectories (hence the optimal trajectories) without error \cite{StuartHumpries96,40}.

\subsection{Problem Formulation and Contributions}

Motivated by the need to understand the dependency of the optimal cost on parameter, constraints, and perturbations, we formulate a hybrid optimal control problem for hybrid inclusions
and reveal key properties about its regularity and existence of solutions. Specifically, 
we consider hybrid systems described by constrained differential and difference inclusions 
as in \cite{hybridbook,220}, which are given by
\begin{equation}
	\HS
	\left\{
	\begin{aligned}
		\dot{x}		&\in F(x)		& x&\in C\\
		x^+				&\in G(x)		& x&\in D,
	\end{aligned}
	\right.
	\label{eq:H}
\end{equation}
The \textit{flow map}~${F:\reals^n\rightrightarrows \reals^n}$ defines the continuous-time evolution~(\textit{flows}) of the state~$x\in\reals^n$ on the \textit{flow set}~${C\subset\dom F}$.
The \textit{jump map}~${G:\reals^n\rightrightarrows \reals^n}$ defines the discrete transitions~(\textit{jumps}) of~$x$ on the \textit{jump set}~$D\subset\dom G$. Informally, a solution of $\HS$ is a function~$(t,j)\mapsto x(t,j)$, where $t$ defines the flow time and $j$ defines the number of jumps. Given a constraint set~${\Omega}\subset\reals^n\times\realsgeq\times\nats\times\reals^n$ and a cost function~$\J:\reals^n\times\realsgeq\times\nats\times\reals^n\to\realsgeq$, the corresponding hybrid optimal control problem we consider is given as follows:
\begin{equation}
\label{eq:optimum}
	\begin{aligned}
		& \underset{x\in\widehat{\sol}_{\HS}}{\text{minimize}}	& & \J{(x(0,0),(T,J),x(T,J))}							\\
		& \text{subject to}														& & (x(0,0),(T,J),x(T,J))\in \Omega,
	\end{aligned}
\end{equation}
where~$\widehat{\sol}_{\HS}$ denotes the set of solutions of~$\HS$ with compact domains (see Section~\ref{sec:background}), and~$(T,J)$ denotes the terminal (hybrid) time of~$x$.\footnote{The notion of solutions is made precise in the next section. For now, we note that solutions are parametrized by hybrid time~$(t,j)$, where~$t\geq 0$ is the ordinary time elapsed and~$j\in\{0,1,\dots\}$ is the number of jumps that has occurred.} {When the cost function~$\J$ depends only on the terminal point~$x(T,J)$ and the constraint set~$\Omega$ is of the form~$\{x_0'\}\times\{(T',J')\}\times X$, this is a standard initial value problem in Mayer form with terminal constraints. Problems similar to~\eqref{eq:optimum} (e.g. variable time with boundary constraints) are considered in~\cite{clarke}; see, for example, Problem (OC1) therein.
}

Our choice of the relatively simplistic structure of optimization problem in~\eqref{eq:optimum} is motivated by the possibility of passing from more general problems to the one in~\eqref{eq:optimum}.\footnote{This also justifies our use of the term optimal ``control'' over alternative descriptors, e.g., calculus of variations.} For example, given the Bolza cost functional in~\cite{acc} for controlled hybrid equations, which includes stage costs for flows and jumps, one can pass to a Mayer cost functional as in~\eqref{eq:optimum} by augmenting the dynamics with an additional state representing the running cost. The continuous/discrete-time analogues of this trick are well known in the literature and can be found in standard references on optimal control, such as~\cite{liberzon}. For the control inputs, we refer to Filippov's lemma (e.g.~\cite[Corollary~23.4]{clarke}), which establishes equivalence between solutions of a controlled differential equation and the corresponding differential inclusion. Finally, we observe that state constraints aside from endpoint constraints are omitted in~\eqref{eq:optimum}, since these can be embedded in the flow set~$C$ and jump set~$D$, as noted in~\cite{9147972}. A similar approach has been taken in~\cite{wolenski}, where the author studies the continuous-time counterpart of~\eqref{eq:optimum} to characterize the value function.

This paper reveals the following key properties of the hybrid optimal control problem in \eqref{eq:optimum}, using recently developed notions of well posedness for hybrid systems (\cite{arxiv}) and their applications to reachable sets:
\begin{enumerate}
\item existence of optimal solutions;
\item upper semicontinuous dependence of the optimal cost on initial conditions, time, magnitude of perturbations, and constraints of the optimal control problem;
\item lower semicontinuous dependence of the optimal cost on initial conditions, time, magnitude of perturbations, and constraints of the optimal control problem; 
\item continuous dependence of the optimal cost on initial conditions, time, magnitude of perturbations, and constraints of the optimal control problem; 
\item outer/upper semicontinuous dependence of optimal solutions on initial conditions, time, magnitude of perturbations, and constraints of the optimal control problem.
\end{enumerate}

The results are illustrated in multiple examples in Section~\ref{sec:ex}. The first two results require the hybrid system in question to have the so-called ``outer well-posedness'' property, which is guaranteed under mild regularity conditions. Lower semicontinuity of the optimal cost requires ``inner well-posedness'', guaranteed under a combination of regularity, tangent cone, and geometric conditions, and also necessitate some assumptions on the structure of solutions. Consequently, combining inner/outer well-posedness properties with this assumption lead to continuity of the optimal cost and upper semicontinuity of optimal solutions. Importantly, a) the aforementioned perturbations include perturbations to the right-hand sides of the differential/difference inclusions defining the hybrid system, as well as the associated constraint sets, and b) as shown in Section~\ref{sec:remarks}, when the assumption on the solution cannot be satisfied, continuous (respectively, outer/upper semicontinuous) approximations of the optimal cost (respectively, solutions) are still possible.

\subsection{Organization of the Paper}

Section~\ref{sec:background} pertains to basic concepts of
 hybrid inclusions and set-valued analysis.
Section~\ref{sec:background-WP} presents an overview about
well-posed hybrid systems.
Section~\ref{sec:MainResults} presents  the main results about
continuity of the optimal cost and upper semicontinuity of optimal solutions.
Section~\ref{sec:remarks} makes remarks about the assumptions involved in the continuity properties established in Section~\ref{sec:MainResults}.
Section~\ref{sec:ex} presents two examples to which the main results
are applied.

\section{Preliminaries}
\label{sec:background}
Throughout the paper,~$\reals$ denotes real numbers,~$\realsgeq$ nonnegative reals, and~$\naturals$ nonnegative integers. The 2-norm is denoted~$|.|$. Given a pair of sets~$S_1,S_2$,~$S_1\subset S_2$ indicates~$S_1$ is a subset of~$S_2$, not necessarily proper. Let~$\A\subset\reals^n$ be nonempty. The distance of a vector~$x\in\reals^n$ to the set~$\A$ is~${\Anorm{x}:=\inf_{a\in\A}|x-a|}$. The closed unit ball in~$\reals^n$ centered at the origin is denoted~$\ball$, $r\ball$ is the closed ball of radius~$r$ centered at the origin, and~$\A+r\ball$ is the set of all~$x$ such that~$|x-a|\leq r$ for some~$a\in\A$. The closure, interior, and boundary of a set~${S\subset\reals^n}$ are denoted~${\closure S}$,~${\interior S}$, and~${\partial S}$. The domain of a set-valued mapping~$M:S\rightrightarrows \reals^m$, denoted~$\dom M$, is the set of all~$x\in S$ such that~$M(x)$ is nonempty. Given a set~$S'\subset S$,~$M|_{S'}$ denotes the restriction of~$M$ to~$S'$.

\subsection{Hybrid Inclusions: Solutions and Reachable Sets}
\label{sec:SolutionsAndReachSets}

We introduce the concept of solution to the hybrid system in~\eqref{eq:H}, whose data is the 4-tuple $(C,F,D,G)$ and, at times, we refer to it using the notation~$\HS=(C,F,D,G)$. Solutions of the hybrid system~$\HS$ belong to a class of functions called \textit{hybrid arcs}. Hybrid arcs are parametrized by hybrid time~$(t,j)$, where~$t\in\realsgeq$ denotes the ordinary time and~$j\in\nats$ denotes the number of jumps. A function~$x$ mapping a subset of~$\realsgeq\times\naturals$ to~$\reals^n$ is a hybrid arc if 1) its domain, denoted~$\dom x$, is a \textit{hybrid time domain}, and 2) it is locally absolutely continuous on each connected component of~$\dom x$. Formally, a set $E \subset\realsgeq\times\naturals$ is a hybrid time domain if for every~$(T,J)\in E$, there exists a nondecreasing sequence~$\{t_j\}_{j=0}^{J+1}$ with~$t_0=0$ such that~$E \cap \left([0,T]\times \{0,1,\dots,J\}\right) = \cup_{j=0}^{J}\left([t_j,t_{j+1}]\times \{j\}\right)$. Then, a function~$x:\dom x\to\reals^n$ is a hybrid arc if~$\dom x$ is a hybrid time domain and for every~$j\geq 0$, the function~$t\mapsto x(t,j)$ is locally absolutely continuous on the interval~$I^j:=\{t:(t,j)\in \dom x\}$. A hybrid arc~$x$ satisfying the dynamics in~\eqref{eq:H} is a solution of the hybrid system~$\HS$ if it satisfies the initial condition constraint~$x(0,0)\in \closure(C)\cup D$~\cite[Definition~2.6]{hybridbook}.

A hybrid arc~$x$ is called \textit{complete} if its domain is unbounded. It is called \textit{bounded} if its range is bounded. It is said to have \textit{finite escape time} if~$x(t,J)$ tends to infinity as~$t$ tends to~$T$ from the left. If the domain of~$x$ is compact, we say that~$(T,J)\in\dom x$ is the \textit{terminal (hybrid) time} of~$x$ if~$t\leq T$ and~$j\leq J$ for all~$(t,j)\in\dom x$. Similarly,~$T$ is referred to as the \textit{terminal ordinary time} of~$x$. The same terminology is used for hybrid arcs that are solutions of the hybrid system~$\HS$; e.g., a solution~$x$ of~$\HS$ is bounded if its range is bounded.

A solution~$x$ of the hybrid system~$\HS$ is \textit{maximal} if it cannot be extended to another solution. The notation~$\sol_{\HS}(S)$ refers to the set of all maximal solutions~$x$ of~$\HS$ originating from~$S$ (i.e.,~${x(0,0)\in S}$ for every~$x\in\sol_{\HS}(S)$), and~$\sol_{\HS}:=\sol_{\HS}(\reals^n)$. If every~$x\in\sol_{\HS}(S)$ is bounded or complete, we say that~$\HS$ is \textit{pre-forward complete from~$S$}. We say that~$t$ is a \textit{jump time} of~$x$ if there exists~$j$ such that~$(t,j),(t,j+1)\in\dom x$. The notation~$\widehat{\sol}_{\HS}$ in \eqref{eq:optimum} denotes the set of all solutions of~$\HS$ (not necessarily maximal) with compact hybrid domains; i.e., $\dom x$ is compact for every~$x\in\widehat{\sol}_{\HS}$. Note that every such~$x$ has a terminal hybrid time~$(T,J)\in\realsgeq\times\nats$.

Given an initial~$x_0$ condition and a hybrid time~$(T,J)$, we define the reachable set of the hybrid system $\HS$
as the set of points reached by solutions originating from~$x_0$ at hybrid time~$(T,J)$.

\begin{defn}[Reachable Set Mappings]
\label{def:reach}
Given a hybrid system~$\HS=(C,F,D,G)$, the reachable set mapping~$\reach_{\HS}:(\closure(C)\cup D)\times\realsgeq\times\nats\rightrightarrows\reals^n$ of~$\HS$ is the set-valued mapping that associates with every~$x_0$,~$T$, and~$J$, the reachable set of~$\HS$ from~$x_0$ at time~$(T,J)$, i.e., $\reach_{\HS}(x_0,T,J):=\{x(T,J): x\in\sol_{\HS}(x_0), (T,J)\in\dom x\}$.
\end{defn}

\subsection{Limits, Semicontinuity, and Boundedness of Set-Valued Maps}

Let~$S\subset\reals^n$,~$x\in \closure S$, and consider a set-valued mapping~$M:S\rightrightarrows \reals^m$. The \textit{inner limit} of~$M$ as~$x'$ tends to~$x$,~$\liminf_{x'\to x}M(x')$, is the set of all~$y$ such that for any sequence~$\{x_i\}_{i=0}^{\infty}\in S$ convergent to~$x$, there exist~$\imath\geq 0$ and a sequence~$\{y_i\}_{i=\imath}^{\infty}$ convergent to~$y$ with~$y_i\in M(x_i)$ for all~$i\geq \imath$. The \textit{outer limit} of~$M$ as~$x'$ tends to~$x$,~$\limsup_{{x'\to x}}M(x')$, is the set of all~$y$ for which there exists a sequence~$\{x_i\}_{i=0}^{\infty}\in S$ convergent to~$x$ and a sequence~$\{y_i\}_{i=0}^{\infty}$ convergent to~$y$ with~$y_i\in M(x_i)$ for all~$i\geq 0$. If the inner and outer limits (as~$x'$ tends to~$x$) are equal, the limit of~$M$ as~$x'$ tends to~$x$, denoted~$\lim_{x'\to x}M(x')$, is defined to be equal to them. Limits of sequences of sets are defined in the same manner. Let~$X\subset S$ and~$x\in \closure X$. Then, the mapping~$M$ is \textit{inner semicontinuous} (respectively, \textit{outer semicontinuous}) at~$x$ relative to~$X$ if the inner (respectively, outer) limit of~$M|_{X}$ as~$x'$ tends to~$x$ contains (respectively, is contained in) $M(x)$. It is \textit{continuous} at~$x$ relative to~$X$ if it is both inner and outer semicontinuous at~$x$ relative to~$X$. In addition,~$M$ is \textit{locally bounded at~$x\in X$ relative to~$X$} if there exists~$\varepsilon>0$ such that the set~$M((x+\varepsilon\ball)\cap X)$ is bounded. When these properties hold for all~$x\in X$, we drop the qualifier ``at~$x$'', and if~$X=S$, we drop the qualifier ``relative to~$X$''. These definitions follow~\cite[Definitions~4.1,~5.4, and~5.14]{rockafellarwets}.\footnote{For locally bounded set-valued maps with closed values, outer semicontinuity is equivalent to \textit{upper semicontinuity}~\cite[Definition~1.4.1]{aubin}, see~\cite[Lemma~5.15]{hybridbook}. Inner semicontinuity coincides with \textit{lower semicontinuity}~\cite[Definition~1.4.2]{aubin}.}


\subsection{Graphical Convergence and Closeness of Hybrid Arcs}

A sequence~$\{x_i\}_{i=0}^{\infty}$ of hybrid arcs is said to be \textit{locally eventually bounded} if for any~$\tau\geq 0$, there exist~$\imath\geq 0$ and a compact set~$K$ such that~$x_i(t,j)\in K$ for every~${i\geq \imath}$ and~$(t,j)\in\dom x_i$ with~$t+j\leq\tau$. It is said to \textit{converge graphically} to a mapping~$M:\realsgeq\times\nats\rightrightarrows\reals^n$, called the \textit{graphical limit} of~$\{x_i\}_{i=0}^{\infty}$,  if the sequence~$\{\gph x_i\}_{i=0}^{\infty}$ converges to~$\gph M$ (in the set convergence sense), where~$\gph $ denotes the graph of a set-valued mapping. Graphical convergence is motivated by the fact that solutions of a hybrid system can have different time domains, which renders the uniform norm an insufficient metric to analyze convergence; see~\cite[Chapter~5]{hybridbook}. In lieu of the uniform norm, we use a concept called~$(\tau,\varepsilon)$-closeness, given in Appendix~\ref{sec:extras}.

\section{Background on Well-Posed Hybrid Systems}
\label{sec:background-WP}

Fundamental in our analysis are the various notions of well posedness for hybrid systems. This section provides a brief overview of these notions, to keep the paper self contained.

\subsection{Nominal Well-Posedness}

Roughly speaking, \textit{nominally outer well-posed}\footnote{This notion has previously been referred to in the literature simply as nominal well-posedness; e.g.~\cite[Definition~6.2]{hybridbook}. The new terminology was introduced in~\cite{cdc2020} to accommodate the then novel notion of nominal inner well-posedness.} hybrid systems are those hybrid systems whose solutions depend outer semicontinuously on initial conditions: for a hybrid system~$\HS$ that is nominally outer well-posed on a set~$S$, the graphical limit~$x$ of a locally eventually bounded graphically convergent sequence~$\{x_i\}_{i=0}^{\infty}$ of solutions is itself a solution. The precise definition is recalled below.
\begin{defn}{\cite[Definition 3.2]{cdc2020}}
\label{def:wellposed}
A hybrid system~$\HS$ is said to be \textit{nominally outer well-posed} on a set\footnote{For all notions of well-posedness, for simplicity, we omit the qualifier ``on $S$'' when~$S=\reals^n$. Also, we say ``at $x_0$'' instead of ``on $S$'' if $S=\{x_0\}$ for some~$x_0$.}~$S\subset\reals^n$ if for every graphically convergent sequence of solutions~$\{x_i\}_{i=0}^{\infty}$ of~$\HS$ satisfying~$\lim_{i\to\infty}x_i(0,0)=:x_0\in S$, the following holds:
	\begin{itemize}
		\item if the sequence~$\{x_i\}_{i=0}^{\infty}$ is locally eventually bounded, then the graphical limit~$x$ is a solution of~$\HS$ originating from~$x_0$;
		\item if the sequence~$\{x_i\}_{i=0}^{\infty}$ is not locally eventually bounded, then there exists~$(T,J)\in\realsgeq\times\nats$ such that~$x=M|_{\dom M\cap([0,T)\times\{0,1,\dots,J\})}$ is a solution of~$\HS$ originating from~$x_0$ that escapes to infinity at time~$(T,J)$, where~$M$ is the graphical limit of~$\{x_i\}_{i=0}^{\infty}$.
	\end{itemize}
\end{defn}

See also \cite[Definition 6.2]{hybridbook} and the discussion below \cite[Lemma 2]{arxiv}. In simple words, this property guarantees that small variations in the initial condition does not lead to large changes in the behavior of solutions. Importantly, nominal well-posedness is implied when the data of the system satisfies mild regularity conditions called the \textit{hybrid basic conditions} (\cite[Assumption~6.5]{hybridbook}), see Theorem~\ref{thrm:nomout} in Appendix~\ref{sec:hbc}.

The natural counterpart to the nominal outer well-posedness is called \textit{nominal inner well-posedness}. For hybrid system~$\HS=(C,F,D,G)$ nominally inner well-posed on~$S$, given a bounded or complete solution~$x$ originating from~$S$ and a sequence of initial conditions~$\{\xi_i\}_{i=0}^{\infty}\in\closure(C)\cup D$ convergent to~$x(0,0)$, one can find a locally bounded sequence of solutions~$\{x_i\}_{i=0}^{\infty}$ graphically convergent to~$x$. This is a constructive property, in the sense that it guarantees that a given solution can be approximated by other solutions with small variations in their initial conditions. Sufficient conditions guaranteeing nominal inner well-posedness are provided in Theorem~\ref{thrm:viab} in Appendix~\ref{sec:hbc}. 
\begin{defn}{\cite[Definition 5]{arxiv}}
\label{def:mydef}
A hybrid system~$\HS=(C,F,D,G)$ is said to be nominally inner well-posed on a set~$S\subset\reals^n$ if for every solution~$x$ of $\HS$ originating from~$S$, the following holds:
\begin{enumerate}[label={($\star$)},leftmargin=*]
	\item \label{item:star} given any sequence~$\{\xi_{i}\}_{i=0}^{\infty}\in\closure(C)\cup D$ convergent to~$x(0,0)$, for every~$i\geq 0$, there exists a solution~$x_i$ of~$\HS$ originating from~$\xi_i$ such that
	\begin{enumerate}[label={(\alph*)},leftmargin=*]
		\item \label{item:stara} if~$x$ is bounded or complete and~$\dom x$ is closed, then the sequence~$\{x_i\}_{i=0}^{\infty}$ is locally eventually bounded and graphically convergent to~$x$;
		\item \label{item:starb} if~$x$ escapes to infinity at hybrid time~$(T,J)$, then the sequence~$\{x_i\}_{i=0}^{\infty}$ is not locally eventually bounded but graphically convergent to a mapping~$M$ such that~$x=M|_{\dom M\cap([0,T)\times\{0,1,\dots,J\})}$.
	\end{enumerate}
\end{enumerate} 
\end{defn}

\subsection{Well-Posedness}

Consider a hybrid system~$\HS_{\delta}=(C_{\delta},F_{\delta},D_{\delta},G_{\delta})$ parametrized by a scalar~$\delta\in(0,1)$. The notions of outer and inner well-posedness are concerned with the behavior of solutions as the parameter~$\delta$ tends to zero. Roughly speaking, given a hybrid system~$\HS$, a family\footnote{For simplicity of notation, we drop the subscript~$\delta\in(0,1)$ when referring to such families.} of hybrid systems~$\{\HS_{\delta}=(C_{\delta},F_{\delta})\}_{\delta\in(0,1)}$ is said to be an inner well-posed perturbation of~$\HS$ if given a bounded or complete solution~$x$ of~$\HS$ originating from~$S$, one can find a locally bounded sequence of solutions~$\{x_i\}_{i=0}^{\infty}$ of this family graphically convergent to~$x$. Just like nominal inner well-posedness, this property guarantees that a given solution of the nominal system can be approximated with small variations in their initial condition, provided the perturbation parameter~$\delta\in(0,1)$ is also small. For sufficient conditions for inner well-posedness, see \cite{arxiv}. 

\begin{defn}[Inner Well-Posed Perturbations]
\label{def:iwppert}
A family of hybrid systems~$\{\HS_{\delta}=(C_{\delta},F_{\delta},D_{\delta},G_{\delta})\}$ is said to be an inner well-posed perturbation of a hybrid system~$\HS$ on a set~$S$ if~$S\cap(\closure(C)\cup D)\subset\liminf_{\delta\to 0}\closure(C_{\delta})\cup D_{\delta}$, and for every solution~$x$ of $\HS$ originating from~$S$, the following hold:
\begin{enumerate}[label={($\diamond$)},leftmargin=*]
	\item \label{item:asterisk} given any sequence~$\{\delta_i\}_{i=0}^{\infty}\in(0,1)$ convergent to zero and any sequence~$\{\xi_{i}\}_{i=0}^{\infty}$ convergent to~$x(0,0)$ with~$\xi_i\in\closure(C_{\delta_i})\cup D_{\delta_i}$ for all~$i\geq 0$, for every~$i\geq 0$, there exists a solution~$x_i$ of~$\HS_{\delta_i}$ originating from~$\xi_i$ such that~\ref{item:stara} and \ref{item:starb} in Definition~\ref{def:mydef} hold.
\end{enumerate}
\end{defn}

Outer well-posedness, being a property tailored towards robustness, considers a specific family of hybrid systems, namely, those given by~$\rho$-perturbations defined in Appendix~\ref{sec:extras}, and requires the analogue of the graphical convergence property for nominal outer well-posedness to hold for \textit{all}~$\rho$-perturbations with continuous function~$\rho$. Every outer well-posed system is nominally outer well-posed, and hybrid basic conditions guarantee outer well-posedness as well as nominal outer well-posedness; see Theorem~\ref{thrm:nomout} in Appendix~\ref{sec:hbc}. The relevant definitions (\cite[Definitions 6.27 and 6.29]{hybridbook}) are recalled below for completeness. 

\begin{defn}[$\rho$-Perturbation]
\label{def:rho}
Given a hybrid system~$\HS=(C,F,D,G)$ and a function~$\rho:\reals^n\to\realsgeq$, the~$\rho$-perturbation of~$\HS$ is the hybrid system~$\HS^{\rho}$ with data~$(C^{\rho},F^{\rho},D^{\rho},G^{\rho})$, where $C^{\rho}=\{x: (x+\rho(x)\ball)\cap C\neq \varnothing\}$, $D^{\rho}=\{x: (x+\rho(x)\ball)\cap D\neq \varnothing\}$, and 
	\begin{equation*}
		\begin{aligned}
			F^{\rho}(x)&=\closure(\con F((x+\rho(x)\ball)\cap C))+\rho(x)\ball,\\
			G^{\rho}(x)&=\{z: z\in y+\rho(y)\ball, y\in G((x+\rho(x)\ball)\cap D)\},
		\end{aligned}
	\end{equation*}
for all~$x\in\reals^n$, where $\con$ denotes the convex hull. Moreover, given any~$\delta\in(0,1)$,~$\HS^{\delta\rho}$ denotes the~$\tilde{\rho}$-perturbation of~$\HS$, where~$\tilde{\rho}$ is the function~$x\mapsto\delta\rho(x)$.
\end{defn}

\begin{defn}[Outer Well-Posedness]
A hybrid system~$\HS$ is said to be \textit{outer well-posed} on a set~$S\subset\reals^n$ if for every continuous function~$\rho$, every positive sequence~$\{\delta_i\}_{i=0}^{\infty}$, and every graphically convergent sequence of hybrid arcs~$\{x_i\}_{i=0}^{\infty}$ such that~$x_i$ is a solution of~$\HS^{\delta_i\rho}$ and~$\lim_{i\to\infty}x_i(0,0)=:x_0\in S$, the following holds:
	\begin{itemize}
		\item if the sequence~$\{x_i\}_{i=0}^{\infty}$ is locally eventually bounded, then the graphical limit~$x$ is a solution of~$\HS$ originating from~$x_0$;
		\item if the sequence~$\{x_i\}_{i=0}^{\infty}$ is not locally eventually bounded, then there exists~$(T,J)\in\realsgeq\times\nats$ such that~$x=M|_{\dom M\cap([0,T)\times\{0,1,\dots,J\})}$ is a solution of~$\HS$ originating from~$x_0$ that escapes to infinity at time~$(T,J)$, where~$M$ is the graphical limit of~$\{x_i\}_{i=0}^{\infty}$.
	\end{itemize}
\end{defn}

The relationship between general families of hybrid systems~$\HS_{\delta}$ and~$\rho$-perturbations are made concrete by the notion of \textit{domination by a $\rho$-perturbation}. Essentially, given a function~$\rho$, a family~$\HS_{\delta}$ is dominated by the $\rho$-perturbation of~$\HS$ if it can be overapproximated by the $\rho$-perturbation.
\begin{defn}[Domination by a $\rho$-Perturbation]
A family of hybrid systems~$\{\HS_{\delta}=(C_{\delta},F_{\delta},D_{\delta},G_{\delta})\}$ is said to be dominated by the $\rho$-perturbation of a hybrid system~$\HS$ if for every~$\delta\in(0,1)$, $C_{\delta}\subset C^{\delta\rho}$, $F_{\delta}(x)\subset F^{\delta\rho}(x)$ for all~$x\in \reals^n$, $D_{\delta}\subset D^{\delta\rho}$, and $G_{\delta}(x)\subset G^{\delta\rho}(x)$ for all~$x\in \reals^n$, where~$(C^{\delta\rho},F^{\delta\rho},D^{\delta\rho},G^{\delta\rho})$ is the data of the $\delta\rho$-perturbation of~$\HS$; see Definition~\ref{def:rho}.
\end{defn}

\section{Key Properties of the Hybrid Optimal Control Problem: Existence and Dependency}
\label{sec:MainResults}

We present our main results on existence of optimal solutions, along with regularity of the optimal cost and the set of optimal solutions to the hybrid optimal control problem in \eqref{eq:optimum}. The proofs reveal that the regularity properties of reachable set mappings are closely related to the aforementioned regularity properties of the optimal costs, since~\eqref{eq:optimum} can equivalently be represented as a finite-dimensional minimization problem over an appropriate reachable set of an augmented hybrid system. Our approach relies on exploiting this link using Berge's maximum theorem~\cite[Theorem~1.4.16]{aubin}. A stronger version of the theorem further enables us to conclude regularity, more precisely, upper/outer semicontinuity of the set of optimal solutions.

Before introducing our results, we note that the terminology concerning~\eqref{eq:optimum} is standard: the hyrid optimal control problem~\eqref{eq:optimum} is said to be \textit{feasible} if there exists a solution~$x$ of~$\HS$ that respects the constraint in~\eqref{eq:optimum}, with~$x$ referred to as a \textit{feasible solution} of~\eqref{eq:optimum}. The \textit{optimal cost} of the problem, denoted $\J^*_{\HS}(\Omega)$, is the infimum of~$\J$ over all feasible solutions, with~$\J^*_{\HS}(\Omega):=\infty$ if~\eqref{eq:optimum} is not feasible, i.e.,
\[
	\J^*_{\HS}(\Omega):=\inf_{\substack{x\in\widehat{\sol}_{\HS}\\ (x(0,0),(T,J),x(T,J))\in \Omega}}\J(x(0,0),(T,J),x(T,J)),
\]
where~$(T,J)$ denotes the terminal time of~$x$. A feasible solution of~\eqref{eq:optimum} that attains the infimum~$\J^*_{\HS}(\Omega)$ is said to be an \textit{optimal solution} of~\eqref{eq:optimum}.

\subsection{Existence of Optimal Solutions and Upper Semicontinuity of the Optimal Cost}
\label{sec:opti-existence}
Within the setting of nominally outer well-posed hybrid systems, it is fairly straightforward to prove existence of optimal solutions under standard regularity conditions. This approach differs from the one in~\cite{goebel2019}, in that it requires no assumptions on the corresponding optimal control problem for the underlying continuous-time system.

\begin{thm}[Existence of Optimal Solutions]
\label{thrm:exist}
Let~$\HS=(C,F,D,G)$ be a hybrid system. Given a compact set~$K$, suppose that~$\HS$ is nominally outer well-posed on~$K$ and pre-forward complete from~$K$. Then, given a closed constraint set~${\Omega}$ and a cost function~$\J$ that is lower semicontinuous on~$\Omega$, there exists an optimal solution of the optimal control problem~\eqref{eq:optimum} if it is feasible and~$\Omega\subset K\times\T\times X$ for a compact set~$\T$.
\end{thm}
\begin{pf}
Since the set~$\Omega\subset K\times\T\times X$ is closed and the set~$ K\times\T$ is compact, the projection of~$\Omega$ onto~$\reals^n\times(\realsgeq\times\nats)$, denoted~$\mathcal{C}$, is compact. Observe that given a feasible solution~$x$, $(x(0,0),(T,J))\in\mathcal{C}$, where~$(T,J)$ is the terminal time of~$x$. Construct an augmented hybrid system~$\HS'$ with state~$z:=(\eta,s,i,x)$, where~$\eta$ represents the initial condition,~$(s,i)$ represents hybrid time, and
~$x$ evolves according to the dynamics of~$\HS$, given by
\begin{equation}
\label{eq:aug}
	\HS'
	\left\{
	\begin{aligned}
		\dot{z}		&\in \{0\}\times\{1\}\times\{0\}\times F(x)		& z&\in C'\\
		z^+				&\in \{\eta\}\times\{s\}\times\{i+1\}\times G(x)		& z&\in D',
	\end{aligned}
	\right.
\end{equation}
with~$C':=\reals^n\times\realsgeq\times\nats\times C$ and~$D':=\reals^n\times\realsgeq\times\nats\times D$. Since~$\HS$ is nominally outer well-posed on~$K$, it is straightforward to show nominal outer well-posedness of~$\HS'$ on the compact set~$K':=\{z:\eta=x\in K, s=i=0\}$. Similarly, since~$\HS$ is and pre-forward complete from~$K$, one can show that~$\HS'$ is pre-forward complete from~$K'$. From these two facts and~\cite[Proposition~4.2]{cdc2020}, the reachable set~$\reach_{\HS'}(\mathcal{C}')$ is compact, where~$\mathcal{C}':=\{(z,T,J):z\in K', (\eta,T,J)\in\mathcal{C}\}$. Consequently, the intersection of $\reach_{\HS'}(\mathcal{C}')$ and $\Omega$ is compact. The optimal control problem~\eqref{eq:optimum} can then be recast as the minimization of the cost function~$\J$ on this intersection. Since the optimal control problem~\eqref{eq:optimum} is feasible, the intersection must be nonempty, and the minimum of~$\J$ on the intersection exists due to lower semicontinuity of~$\J$. Hence, there exists an optimal solution. \hfill $\blacksquare$
\end{pf}


Under the conditions of Theorem~\ref{thrm:exist}, it is also possible to show that the optimal cost depends upper semicontinuously on constraints. This result can be used to show that the \textit{value function}\footnote{The value function corresponding to~\eqref{eq:optimum} is the mapping~$x_0\mapsto\J^*_{\HS}(\Omega)$ in the specific case of~$\Omega=\{x_0\}\times \mathcal{C}$ for some~$\mathcal{C}\subset\realsgeq\times\nats\times\reals^n$.} is upper semicontinuous. As a gentle reminder, upper/lower semicontinuity of (extended) real-valued functions should not be confused with upper/lower semicontinuity of set-valued maps.\footnote{Note that, we use~$\liminf$ to denote the inner limit of sets and set-valued maps, as well as the limit inferior of functions.}

\begin{thm}
\label{thrm:usc}
Let~$\HS$ be a hybrid system and given a compact set~$K$, suppose that~$\HS$ is nominally outer well-posed on~$K$ and pre-forward complete from~$K$. Consider a closed constraint set~${\Omega}$ such that~$\Omega\subset K\times\T\times X$ for a compact set~$\T\subset\realsgeq\times\nats$ and a cost function~$\J$ that is lower semicontinuous on~$\Omega$. Let~$S\subset\reals^m$ be a set containing the origin and~$M:S\rightrightarrows\reals^n\times\realsgeq\times\nats\times\reals^n$ be a set-valued mapping that is locally bounded and outer semicontinuous at the origin, with~$M(0)=\Omega$. Then, the function~$\theta\mapsto\J^*_{\HS}(M(\theta))$ is upper semicontinuous at the origin.
\end{thm}
\begin{pf}
The proof follows similar ideas as the proof of Theorem~\ref{thrm:exist}. Local boundedness and outer semicontinuity of~$M$ at the origin implies that given the projection~$\mathcal{C}(\theta)$ of~$M(\theta)$ onto~$\reals^n\times(\realsgeq\times\nats)$, the mapping ~$\mathcal{C}$ is locally bounded and outer semicontinuous at the origin. Let~$K'(\theta):=\{z=(\eta,s,i,x):\eta=x\in \Pi(M(\theta)), s=i=0\}$, where~$\Pi$ is the canonical projection onto the first coordinate, which satisfies~$\Pi(M(0))\subset K$, and note that~$K'$ is locally bounded and outer semicontinuous at the origin (due to~\cite[Proposition~5.52]{rockafellarwets}). Let $\mathcal{C}'(\theta):=\{(z,T,J):z\in K'(\theta), (\eta,T,J)\in\mathcal{C}(\theta)\}$, which is also locally bounded and outer semicontinuous at the origin. Now, construct the augmented hybrid system~$\HS'$ in~\eqref{eq:aug}, which is nominally outer well-posed on the set~$K'(0)$, and note that the reachable set mapping~$\reach_{\HS'}$ is outer semicontinuous and locally bounded at~$\mathcal{C}'(0)$ by~\cite[Theorem~4.1]{cdc2020}. Then, the mapping from~$\theta$ to~$\reach_{\HS'}(\mathcal{C}'(\theta))\cap M(\theta)$ is also outer semicontinuous and locally bounded at the origin. Equivalently, it is upper semicontinuous~\cite[Lemma~5.15]{hybridbook}. Recasting the optimal control problem as the minimization of~$\J$ on this intersection (for each~$\theta)$, Berge's maximum theorem~\cite[Theorem~1.4.16]{aubin} is applicable, and lower semicontinuity of~$\J$, combined with the upper semicontinuity of the aforementioned intersection, leads to upper semicontinuity of~$\theta\mapsto\J^*_{\HS}(M(\theta))$.\hfill $\blacksquare$
\end{pf}


For a fixed-time initial value problem without terminal constraints (i.e., the set $\Omega$ in \eqref{eq:optimum} is of the form~$\{x_0\}\times\{(T,J)\}\times\reals^n$), one can simply take~$M(x'_0\,T',J')=\{x'_0\}\times\{(T',J')\}\times S$ and invoke Theorem~\ref{thrm:usc} to conclude upper semicontinuity of the value function, where~$S$ is an arbitrary compact set that contains the reachable set~$\reach_{\HS}(x_0,T,J)$.\footnote{The reachable set is bounded, in fact compact; see~\cite[Proposition~4.2]{cdc2020}.} Moreover, upper semicontinuous dependence of the value function on initial conditions can easily be extended to show upper semicontinuous dependence on the magnitude of perturbations on the hybrid system and the terminal constraint.

\begin{thm}
\label{thrm:youesssea}
Let~$\HS$ be a hybrid system and given a compact set of initial conditions~$K$, suppose that~$\HS$ is outer well-posed on~$K$ and pre-forward complete from~$K$. Consider a closed constraint set~${\Omega}$ satisfying~$\Omega\subset K\times\T\times\reals^n$ for a compact set~$\T\subset\realsgeq\times\nats$. Let~$\{\HS_{\delta}\}$ be a family of hybrid systems dominated by the~$\rho$-perturbation of~$\HS$ for some continuous function~$\rho$. Moreover, let~$S\subset\reals^m$ be a set containing the origin and~$M:S\rightrightarrows\reals^n\times\realsgeq\times\nats\times\reals^n$ be a set-valued mapping that is locally bounded and outer semicontinuous at the origin, with~$M(0)=\Omega$, and for each~$\theta\in S$, let~$\J(.;\theta)$ be a cost function. Suppose that the function~$(\xi,\theta)\mapsto\J(\xi;\theta)$ is lower semicontinuous at~$(\xi',0)$ for all~$\xi'\in\Omega$. Then, the function~$(\delta,\theta)\mapsto\J^*_{\HS_{\delta}}(M(\theta);\theta)$\footnote{By this, we denote the value function of an optimal control problem with hybrid system~$\HS_{\delta}$, constraint set~$M(\theta)$, and cost function~$\xi\mapsto\J(\xi;\theta)$.} with~$\J^*_{\HS_{0}}(M(\theta);\theta):=\J^*_{\HS}(M(\theta);\theta)$ for all~$\theta\in\reals^m$, is upper semicontinuous at the origin.
\end{thm}

The proof of Theorem~\ref{thrm:youesssea} is almost the same as that of Theorem~\ref{thrm:usc}. It is omitted for brevity. The required outer semicontinuity and local boundedness properties for the reachable set is proved in \cite[Theorem 35]{arxiv}.

\subsection{Continuity of the Optimal Cost and Outer/Upper Semicontinuity of Optimal Solutions}

\label{sec:opti}
Lower semicontinuous, and more strongly, continuous dependence of the optimal cost on constraints and perturbations can similarly be established within the setting of inner well-posedness. Remarkably, the assumptions we use to prove these properties elegantly lead to outer/upper semicontinuous dependence of the set of optimal solutions on constraints and perturbations as well; see Theorem~\ref{thrm:coupdegrace}. In establishing these stronger results, for simplicity and brevity, we focus on fixed-time initial value problems without terminal constraints. That is, we develop our results by focusing on~\eqref{eq:optimum} with the constraint set~$\Omega=\{x_0\}\times\{T\}\times\{J\}\times\reals^n$ (fixed initial value, fixed time, no terminal constraints).

The main results in this subsection consider perturbations to both the hybrid system and the cost function, i.e., they consider parametrized families of hybrid systems and cost functions. Results concerning the nominal case are recovered as immediate corollaries. One key assumption that we make is that the points belonging to the reachable set~$\reach_{\HS}(x_0,T,J)$ correspond to maximal solutions of the hybrid system~$\HS$ that originate from~$x_0$ and do not jump or terminate at ordinary time~$T$. This allows us to conclude lower semicontinuous, and where appropriate, continuous dependence of the reachable set mappings on their arguments and parameters. As shown below, for lower semicontinuity of the optimal cost, it suffices to assume inner well-posedness.

\begin{thm}
\label{thrm:yolo}
Let~$\HS$ be a hybrid system. Given an initial condition~$x_0$ and~$(T,J)\in\realsgeq\times\nats$, suppose that the reachable set~$\reach_{\HS}(x_0,T,J)$ is nonempty, and for every~$\xi\in\reach_{\HS}(x_0,T,J)$, there exists~$x\in\sol_{\HS}(x_0)$ such that~$\xi=x(T,J)$ and~$T$ is not a jump time or the terminal ordinary time of~$x$. Let~$\{\HS_{\delta}=(C_{\delta}, F_{\delta}, D_{\delta}, G_{\delta})\}$ be an inner well-posed perturbation of~$\HS$ at~$x_0$. Given a set~$S\subset\reals^m$ containing the origin, for each~$\theta\in S$, consider a cost function~$\J(.;\theta)$, and suppose that the function~$(\xi,\theta)\mapsto\J(\xi;\theta)$ is upper semicontinuous at~$(\xi',0)$ for all~$\xi'\in\Omega:=\{x_0\}\times\{(T,J)\}\times\reals^n$. For each~$\theta\in S$, let
\begin{equation}
	h_\delta(x_0',T',J';\theta):=\J^*_{\HS_{\delta}}\left(\{x_0'\}\times\{T',J'\}\times\reals^n;\theta\right)
\label{eq:eych}
\end{equation}
for all~$\delta> 0$,~$x_0'\in\closure(C_{\delta})\cup D_{\delta}$,~$(T',J')\in\realsgeq\times\nats$, and~$\theta\in\reals^m$. Then,
\begin{multline}
	\J^*_{\HS}(\Omega;0)=:h(x_0,T,J)\\
	\leq \liminf_{\substack{\delta\to 0,\, \theta\to 0\\(x_0',T',J')\to (x_0,T,J)}}h_\delta(x_0',T',J';\theta).
\label{eq:heyo}
\end{multline}
\end{thm}
\begin{pf}
We go through a reachability analysis as in Theorems~\ref{thrm:exist} and~\ref{thrm:usc}. However, augmentation of the system is not necessary since there are no terminal constraints and the constraints are not mixed. That is, we are interested in instances of problem~\eqref{eq:optimum}, for which the constraint~$(x(0,0),(T,J),x(T,J))\in \Omega$ can be rewritten in the form~$x(0,0)=\xi,(T,J)=(s,i)$. Consequently, given~$\delta$, $\theta$, and $(x_0',T',J')$, the scalar~$h_{\delta}(x_0',T',J;\theta)$ is the minimum of~$\J(.;\theta)$ on the reachable set~$\reach_{\HS_{\delta}}(x_0',T',J')$. Noting that the family of reachable set mappings depend lower semicontinuous on its arguments and the parameter~$\delta$ by \cite[Theorem 36]{arxiv}., i.e.,
\[
\reach_{\HS}(x_0,T,J)\subset\liminf_{\substack{\delta\to 0\\ (x_0',T',J') \to (x_0,T,J)}}\reach_{\HS_\delta}(x_0',T',J'),
\]
and the reachable set~$\reach_{\HS}(x_0,T,J)$ is nonempty, it follows that upper semicontinuity of~$\J$ implies~\eqref{eq:heyo}, c.f. Berge's maximum theorem~\cite[Theorem~1.4.16]{aubin}.\hfill $\blacksquare$
\end{pf}


\begin{cor}
\label{cor:yolo0}
Let~$\HS=(C,F,D,G)$ be a hybrid system. Given an initial condition~$x_0$ and~$(T,J)\in\realsgeq\times\nats$, suppose that the reachable set~$\reach_{\HS}(x_0,T,J)$ is nonempty and~$\HS$ is nominally inner well-posed at~$x_0$. Moreover, suppose that for every~$\xi\in\reach_{\HS}(x_0,T,J)$, there exists~$x\in\sol_{\HS}(x_0)$ such that~$\xi=x(T,J)$ and~$T$ is not a jump time or the terminal ordinary time of~$x$. Consider a cost function~$\J$ that is upper semicontinuous at every~$\xi'\in\Omega:=\{x_0\}\times\{(T,J)\}\times\reals^n$. Let
\begin{equation}
	h(x_0',T',J'):=\J^*_{\HS}\left(\{x_0'\}\times\{T',J'\}\times\reals^n\right)
\label{eq:eych0}
\end{equation}
for all~$x_0'\in\closure(C)\cup D$,~$(T',J')\in\realsgeq\times\nats$. Then,~$h$ is lower semicontinuous at~$(x_0,T,J)$.
\end{cor}

With the additional assumption of outer well-posedness, combining Theorems~\ref{thrm:youesssea} and~\ref{thrm:yolo}, we reach Theorem~\ref{thrm:gme}, which guarantees that the optimal cost can be continuously approximated. Similar to our prior comment regarding Theorem~\ref{thrm:usc}, to be able to invoke Theorem~\ref{thrm:youesssea}, one can consider a ``terminal constraint set'' $S$ containing the reachable set~$\reach_{\HS}(x_0,T,J)$ in its interior. Its immediate corollary, Corollary~\ref{cor:gme}, can also be concluded by combining Theorem~\ref{thrm:usc} and Corollary~\ref{cor:yolo0}.

\begin{thm}[Continuity of the Optimal Cost]
\label{thrm:gme}
Let~$\HS=(C,F,D,G)$ be a hybrid system, and given an initial condition~$x_0$, suppose that~$\HS$ is outer well-posed at~$x_0$ and pre-forward complete from~$x_0$. Let~$\{\HS_{\delta}=(C_{\delta},F_{\delta},D_{\delta},G_{\delta})\}$ be an inner well-posed perturbation of~$\HS$ at~$x_0$ that is dominated by a~$\rho$-perturbation of~$\HS$ for some continuous function~$\rho$. Moreover, given~$(T,J)\in\realsgeq\times\nats$, suppose that the reachable set~$\reach_{\HS}(x_0,T,J)$ is nonempty and~for every~$\xi\in\reach_{\HS}(x_0,T,J)$, there exists~$x\in\sol_{\HS}(x_0)$ such that~$\xi=x(T,J)$ and~$T$ is not a jump time or the terminal ordinary time of~$x$. Given a set~$S\subset\reals^m$ containing the origin, for each~$\theta\in S$, consider a cost function~$\J(.;\theta)$, and suppose that the function~$(\xi,\theta)\mapsto\J(\xi;\theta)$ is continuous at~$(\xi',0)$ for all~$\xi'\in\Omega:=\{x_0\}\times\{(T,J)\}\times\reals^n$. Then, the scalar~$h(x_0,T,J)$ in~\eqref{eq:heyo} and the family of functions~$\{h_{\delta}\}$ in~\eqref{eq:eych} satisfy
\begin{equation}
h(x_0,T,J)= \lim_{\substack{\delta\to 0,\,\theta\to 0\\(x_0',T',J')\to (x_0,T,J)}}h_\delta(x_0',T',J';\theta).
\label{eq:limo}
\end{equation}
\end{thm}

\begin{cor}[Continuity of the Optimal Cost]
\label{cor:gme}
Let~$\HS=(C,F,D,G)$ be a hybrid system, and given an initial condition~$x_0$, suppose that~$\HS$ is outer and inner well-posed at~$x_0$ and pre-forward complete from~$x_0$. Given~$(T,J)\in\realsgeq\times\nats$, suppose that the reachable set~$\reach_{\HS}(x_0,T,J)$ is nonempty and~for every~$\xi\in\reach_{\HS}(x_0,T,J)$, there exists~$x\in\sol_{\HS}(x_0)$ such that~$\xi=x(T,J)$ and~$T$ is not a jump time or the terminal ordinary time of~$x$. Consider a cost function~$\J$ that is continuous at all~$\xi'\in\Omega:=\{x_0\}\times\{(T,J)\}\times\reals^n$. Then, the function~$h$ in~\eqref{eq:eych0} is continuous at~$(x_0,T,J)$.
\end{cor}

Moreover, under the conditions of Theorem~\ref{thrm:gme}, the set of optimal solutions depend on constraints and perturbations in an outer/upper semicontinuous manner, as shown below. In Theorem~\ref{thrm:coupdegrace} below, for fixed~$\delta>0$ and~$\theta\in S$, $\mathcal{O}_{\delta}(x_0',T',J';\theta)$ denotes the set of optimal solutions of the optimal control problem with hybrid system~$\HS_{\delta}$, constraint set~$\{x_0'\}\times\{(T',J')\}\times\reals^n$, and cost function~$\J(.;\theta)$. In the same fashion, $\mathcal{O}(x_0,T,J)$ denotes the set of optimal solutions of the optimal control problem with hybrid system~$\HS$, constraint set~$\{x_0\}\times\{(T,J)\}\times\reals^n$, and cost function~$\J(.;0)$.

\begin{thm}[Optimal Solutions]
\label{thrm:coupdegrace}
Under the conditions of Theorem~\ref{thrm:gme}, the following statements are true.
\begin{description}
	\item[Local Boundedness:] There exist~$\varepsilon>0$ and a compact set~$K$ such that
	\begin{multline}
	\label{eq:leb}
		\delta\in(0,\varepsilon],\,\theta\in\varepsilon\ball,\,\\ x_0'\in x_0+\varepsilon\ball,\,(T',J')\in(T,J)+\varepsilon\ball\\
		\implies x(t,j)\in K \quad
	\end{multline}
	for all $(t,j)\in\dom x$ and $x\in\mathcal{O}_{\HS_{\delta}}(x_0',T',J';\theta)$.
	\item[Outer Semicontinuity:] Let $\{\delta_i\}_{i=0}^{\infty}$ be a positive sequence convergent to zero, $\{\theta_i\}_{i=0}^{\infty}\in S$ be a sequence convergent to zero, and~$\{x'_i\}_{i=0}^{\infty}$ be a graphically convergent sequence of optimal solutions such that $x'_i\in\mathcal{O}_{\HS_{\delta_i}}(\xi_i,T_i,J_i;\theta_i)$ for all~$i\geq 0$. Then, if the sequences~$\{\xi_i\}_{i=0}^{\infty}$ and~$\{(T_i,J_i)\}_{i=0}^{\infty}$ converge to~$x_0$ and~$(T,J)$, respectively,~$x\in\mathcal{O}_{\HS}(x_0,T,J)$,	where~$x$ is the graphical limit of~$\{x'_i\}_{i=0}^{\infty}$.
	\item[Upper Semicontinuity:] For all~$\tau\geq 0$ and~$\varepsilon>0$, there exists~$\eta>0$ such that the following holds: for every~$\delta\in(0,\eta]$, $\theta\in \eta\ball$, $x'_0\in x_0+\eta\ball$, $(T',J')\in(T,J)+\eta\ball$, and~$x'\in\mathcal{O}_{\HS_{\delta}}(x'_0,T',J';\theta)$, there exists~$x\in\mathcal{O}_{\HS}(x_0,T,J)$ such that~$x$ and~$x'$ are $(\tau,\varepsilon)$-close.
\end{description}
\end{thm}
\begin{pf}
Local boundedness of the optimal solutions in the sense of \eqref{eq:leb} is a direct result of outer well-posedness of~$\HS$ and the fact that the family of hybrid systems~$\{\HS_{\delta}\}$ are dominated by a~$\rho$-perturbation of~$\HS$ for a continuous function~$\rho$. It can be concluded by upper semicontinuous dependence of solutions on initial conditions and perturbations (\cite[Proposition~6.34]{hybridbook}) and compactness of the reachable set of~$\HS$ over compact hybrid time horizons~\cite[Proposition~4.2]{cdc2020}. Alternatively, one can invoke Theorem 35 or 37 in \cite{arxiv}.

The second statement can be interpreted as ``well-posedness'' of the optimal solutions. To prove this statement, we first note that by well-posedness and~\cite[Lemma~2]{arxiv}, the graphical limit~$x$ is a solution of~$\HS$ originating from~$x_0$ with terminal hybrid time~$(T,J)$ and terminal point~$x(T,J)=\lim_{i\to\infty}x'_i(T_i,J_i)$. To show optimality, we go through the reachability analysis discussed in the proof of Theorem~\ref{thrm:yolo} and recall that given~$\delta>0$, $\theta\in S$, and~$(x_0',T',J')$, the optimal control problem is equivalent to minimizing the cost~$\J(.;\theta)$ on the reachable set~$\reach_{\HS_{\delta}}(x'_0,T',J')$. For every~$\delta> 0$, $\theta\in S$, and~$(x_0',T',J')$, let 
\[
	M_{\delta}^{\theta}(x'_0,T',J'):=\argmin_{\xi\in\reach_{\HS_{\delta}}(x'_0,T',J')} \J(\xi;\theta),
\]
and similarly, for every~$(x_0',T',J')$, let
\[
M(x'_0,T',J'):=\argmin_{\xi\in\reach_{\HS}(x'_0,T',J')} \J(\xi;0),
\]
Invoking a stronger version of Berge's maximum theorem (c.f.~\cite[Theorem~5.4.3.]{polak}), the mapping~$M_{\delta}^{\theta}$, which collects the terminal points of optimal solutions, is outer semicontinuous at~$(0,0,x_0,T,J)$, in the sense that
\[
	\limsup_{\substack{\delta\to 0,\,\theta\to 0\\ (x_0',T',J') \to (x_0,T,J) }}M_{\delta}^{\theta}(x_0',T',J')\subset M(x_0,T,J),
\]
since the reachable set mapping is compact-valued and continuous by \cite[Theorem 37]{arxiv}., and the cost function is continuous in both arguments by assumption. Thus,
\[
x(T,J)=\lim_{i\to\infty}x'_i(T_i,J_i)\in M(0,x_0,T,J),
\]
and therefore~$x$ is an optimal solution;~$x\in\mathcal{O}_{\HS}(x_0,T,J)$.

The last statement is proven by contradiction. Assuming that the statement is false, there exists a sequence of optimal solutions~$\{x'_i\}_{i=1}^{\infty}$ such that for every~$i\geq 1$, the following holds: 1) $x'_i\in\mathcal{O}_{\HS_{\delta_i}}(\xi_i,T_i,J_i; \theta_i)$ for some~$\delta_i \leq 1/i$, $\theta\in(1/i)\ball$,~$\xi_i\in x_0+(1/i)\ball$, and~$(T_i,J_i)\in(T,J)+(1/i)\ball$, and 2) no~$x\in\mathcal{O}_{\HS}(x_0,T,J)$ is such that~$x$ and~$x'_i$ are $(\tau,\varepsilon)$-close. As shown earlier, optimal solutions are locally bounded, hence the sequence~$\{x'_i\}_{i=1}^{\infty}$ is locally eventually bounded. Using~\cite[Theorem~6.1]{hybridbook} and without relabeling, we pass to a graphically convergent subsequence. The limit of the sequence, say~$x^*$, is then optimal by our prior conclusion. That is,~$x^*\in\mathcal{O}_{\HS}(x_0,T,J)$. However, by \cite[Theorem~5.25]{hybridbook}, for large enough~$i$, $x^*$ and~$x'_i$ are $(\tau,\varepsilon)$-close, which is a contradiction.\hfill $\blacksquare$
\end{pf}


Similarly, in the following, $\mathcal{O}(x_0',T',J')$ denotes the set of optimal solutions of the optimal control problem with hybrid system~$\HS$, constraint set~$\{x_0'\}\times\{(T',J')\}\times\reals^n$, and cost function~$\J$.
\begin{cor}[Optimal Solutions]
\label{cor:coupdegrace}
Under the conditions of Corollary~\ref{cor:gme}, the following statements are true.
\begin{description}
	\item[Local Boundedness:] There exists~$\varepsilon>0$ and a compact set~$K$ such that
	\[
		x_0'\in x_0+\varepsilon\ball,\,(T',J')\in(T,J)+\varepsilon\ball
		\implies x(t,j)\in K \quad
	\]
	for all $(t,j)\in\dom x$ and $x\in\mathcal{O}_{\HS}(x_0',T',J')$.
	\item[Outer Semicontinuity:] Let~$\{x'_i\}_{i=0}^{\infty}$ be a graphically convergent sequence of optimal solutions such that $x'_i\in\mathcal{O}_{\HS}(\xi_i,T_i,J_i)$ for all~$i\geq 0$. Then, if the sequences~$\{\xi_i\}_{i=0}^{\infty}$ and~$\{(T_i,J_i)\}_{i=0}^{\infty}$ converge to~$x_0$ and~$(T,J)$, respectively,~$x\in\mathcal{O}_{\HS}(x_0,T,J)$,	where~$x$ is the graphical limit of~$\{x'_i\}_{i=0}^{\infty}$.
	\item[Upper Semicontinuity:] For all~$\tau\geq 0$ and~$\varepsilon>0$, there exists~$\eta>0$ such that the following holds: for every $x'_0\in x_0+\eta\ball$, $(T',J')\in(T,J)+\eta\ball$, and~$x'\in\mathcal{O}_{\HS}(x'_0,T',J')$, there exists~$x\in\mathcal{O}_{\HS}(x_0,T,J)$ such that~$x$ and~$x'$ are $(\tau,\varepsilon)$-close.
\end{description}
\end{cor}


\section{Remarks on Continuous Approximation of the Optimal Cost and Upper Semicontinuous Approximation of Optimal Solutions}
\label{sec:remarks}
As seen so far, the regularity properties of reachable set mappings are closely related to those of the optimal control problem, since the optimal control problem in~\eqref{eq:optimum} is equivalent to a finite-dimensional minimization problem over an appropriate reachable set. The downside to this approach is that, since reachable set mappings are not continuous in general~\cite{arxiv,semantics}, it might be difficult to argue continuous or lower semicontinuous dependence of the optimal cost with respect to the constraints of the optimal control problem. Indeed, vaguely speaking, the results in Section~\ref{sec:opti} establishing continuity of the optimal cost and semicontinuity properties of optimal solutions require that the parameter $T$ therein is not a jump time or the terminal ordinary time of a solution; see, e.g., the second sentence in the statement of Theorem~\ref{thrm:yolo}. Aside from requiring partial knowledge of solutions, these results cannot guarantee continuity properties of the optimal control problem globally (for example continuity of the optimal cost may not be achievable at certain combinations of initial condition and hybrid time). However, it is easy to show that if the reachable set mappings corresponding to a parametric optimal control problem vary continuously at points related to the problem parameters, the optimal cost varies continuously when there are no terminal constraints (i.e., the constraint set~$\Omega$ in~\eqref{eq:optimum} is of the form~$\mathcal{C}\times\reals^n$). Fortunately, the results reported in \cite{arxiv} show that when the solutions to hybrid systems depend upper and lower semicontinuously on initial conditions and perturbations, the reachable sets can be varied continuously, provided certation class-$\mathcal{K}_{\infty}$ bounds are respected, see Theorems 38 and 40 therein.

For completeness, we include~\cite[Theorem 38]{arxiv} below, restated to omit perturbations to~$\HS$ and only consider a single initial condition and hybrid time, and consequently, state continuity of the optimal cost of the corresponding problem without proof.

\begin{thm}
\label{thrm:cont1}
Let~$\HS$ be a hybrid system, and given an initial condition~$x_0$, suppose that~$\HS$ is nominally inner and outer well-posed at~$x_0$ and pre-forward complete from~$x_0$. Then, for any hybrid time~$(T,J)\in\realsgeq\times\nats$, there exists a class-$\K$ function~$\alpha$ such that for every~$x_0\in K$ and~$(T,J)\in \T$,
\begin{multline*}
	\lim_{\substack{\varepsilon\to 0,\, x_0' \to x_0\\ \varepsilon>0\\ x_0'\in (x_0+\alpha(\varepsilon)\ball)\cap(\closure(C))\cup D)}}\reach_{\HS}(x_0',[\max\{0,T-\varepsilon\},T+\varepsilon],J)\\
	=\reach_{\HS}(x_0,T,J).
\label{eq:wee}
\end{multline*}
\end{thm}

\begin{thm}
Let~$\HS$ be a hybrid system, and given an initial condition~$x_0$, suppose that~$\HS$ is nominally inner and outer well-posed at~$x_0$ and pre-forward complete from~$x_0$. Given a hybrid time~$(T,J)\in\realsgeq\times\nats$, suppose that the reachable set~$\reach_{\HS}(x_0,T,J)$ is nonempty. Consider a cost function~$\J$ that is continuous at all~$\xi'\in\Omega:=\{x_0\}\times\{(T,J)\}\times\reals^n$, and let
\begin{multline*}
h(x_0',\varepsilon)=\J_{\HS}^*(x_0'\times[\max\{0,T-\varepsilon\},T+\varepsilon]\times\{J\}\times\reals^n)\\
 \forall x_0'\in\closure(C)\cup D, \varepsilon\geq 0.
\end{multline*}
Then,
\[
	\lim_{\substack{\varepsilon\to 0,\, x_0' \to x_0\\ \varepsilon>0\\ x_0'\in (x_0+\alpha(\varepsilon)\ball)\cap(\closure(C))\cup D)}}h(x_0',\varepsilon)=h(x_0,0).
\]
\end{thm}


\section{Examples}
\label{sec:ex}
In this section, we consider concrete finite horizon optimization problems for hybrid plants given by
\begin{equation}
	\HS_P
	\left\{
	\begin{aligned}
		\dot{x}_P		&\in F_P(x_P,u)		& (x_P,u)&\in C_P\\
		x_P^+				&\in G_P(x_P,u)		& (x_P,u)&\in D_P,
	\end{aligned}
	\right.
	\label{eq:Hp}
\end{equation}
where $C_P$ is the flow set,
$F_P$ is the flow map,
$D_P$ is the jump set,
and 
$G_P$ is the jump map.
A solution of $\HS_P$ is defined by a pair (called a \textit{solution pair}) $(t,j) \mapsto (x_P(t,j),u(t,j))$
on a hybrid time domain~$\dom(x_P,u)$ satisfying the dynamics of $\HS_P$, in a similar manner as the way 
a solution of the (closed-loop) hybrid system
$\HS$ in \eqref{eq:H} is defined in
Section~\ref{sec:SolutionsAndReachSets}.
Given a solution pair~$(x_P,u)$ with compact domain, the associated cost
 is defined by
\begin{multline}
\label{eq:bestcost}
\hspace{-0.2in}	\left(\sum_{j=0}^{J}\int_{t_j}^{t_{j+1}}L_{C_P}(x_P(t,j),u(t,j))\,dt\right) +  \\
	\left(\sum_{j=0}^{J-1}L_{D_P}(x_P(t_{j+1},j),u(t_{j+1},j))\right)+ V(x_P(T,J)),
\end{multline}
where~$t_j$ is the $j$-th jump time and $(T,J) \in \dom (x_P,u)$ is the terminal time, i.e.,
\[
\dom(x_P,u) = \cup_{j=0}^J\left([t_j,t_{j+1}]\times\{j\}\right)
\]
and~$T=T_{J+1}$. In \eqref{eq:bestcost}, the first term $L_{C_P}$ is the stage cost capturing the cost over intervals of flows, $L_{D_P}$ is the stage cost capturing the cost to jump, and $V$ is the terminal cost. 

The constructions presented above lead to the following finite horizon hybrid optimization problem.

\begin{problem}
\label{prob:opti}
Given a hybrid system~$\HS_P$ as in~\eqref{eq:Hp}, 
a stage cost  for flows $L_{C_P}$, 
 a stage cost for jumps $L_{D_P}$,
a terminal cost $V$, a closed set~$X_P$, a hybrid time~$(T,J)\in\realsgeq\times\nats$,
and an initial condition~$\xi$, 
find a solution pair~$(x_P,u)$ minimizing~\eqref{eq:bestcost} subject to
\begin{itemize}
\item the initial condition constraint $x_P(0,0)=\xi$, and,
\item the terminal constraint $x_P(T,J)\in X_P$.
\end{itemize}
\end{problem}

Note that the flow and jump sets of $\HS_P$ impose constraints that the solution pair needs
to satisfy during flows and jumps, respectively.  In fact, for the solution pair to exist up to hybrid time 
$(T,J)$ it has to belong to $C_P$ and $D_P$: as \cite[Definition 2.29]{220} indicates, $(x_P,u)$ is a solution of $\HS_P$ if
\begin{itemize}
\item
$(x_P(0,0),u(0,0))\in \closure(C_P)\cup D_P$,
\item
For each $j\geq 0$, $(x_P(t,j),u(t,j)) \in C_P$ for all $t \in \interior I^{j}$ and $\dot{x}_P(t,j) \in F_P(x_P(t,j),u(t,j))$ for almost all $t \in I^{j}$, where $I^j:=\{t: (t,j)\in\dom(x_P,u)\}$;
\item 
For each $(t,j)\in \dom (x_P,u)$ such that $(t,j+1)\in \dom (x_P,u)$, $(x_P(t,j),u(t,j))\in D_P$ and
\[
x_P(t,j+1) \in G_P(x_P(t,j),u(t,j)).
\]
\end{itemize}  

Given~$\xi\in\closure(C_P)\cup D_P$ and~$(T,J)\in\realsgeq\times\nats$, $h(\xi,T,J)$ denotes the value of~\eqref{eq:bestcost} given a minimizing solution pair~$(x_P,u)$ of~$\HS_P$ subject to the constraints~$x_P(0,0)=\xi$ and~$x_P(T,J)\in X_P$.

\subsection{Thermostat}

A model capturing the evolution of the temperature of a room controlled by a heater that can either be {\em on} or {\em off} is given by
$$
\dot z = - z + z_o + z_\Delta q,
$$
where $z \in \reals$ is the temperature of the room, 
$z_o$ denotes the effective temperature outside of the room,
$z_\Delta$ represents the capacity of the heater,
and the state $q \in \{0,1\}$ represents whether the heater is {\em on} or {\em off}.  
The value $q = 1$ corresponds to the heater being 
{\em on} and the value $q=0$ indicates that the heater is {\em off}.
Using this continuous-time model, 
we are interested in designing a control algorithm that fulfills the following specifications:
\begin{enumerate}[label=\alph*)]
\item \label{item:Thermostat-Steer} 
steer the temperature to a desired temperature range $[z_{\min},z_{\max}]$, where $z_{\min} < z_{\max}$; and
\item \label{item:Thermostat-Switch} 
minimize the number of {\em on}/{\em off} switches of the heater.
\end{enumerate}
To meet these specifications, we properly define the elements 
in Problem~\ref{prob:opti} and solve it numerically.
The desired steering property can be guaranteed by
selecting the flow cost $L_{C_P}$ as a smooth indicator of the set
$[z_{\min},z_{\max}]$.
The jump cost $L_{D_P}$ can be used to penalize switches from {\em on} to {\em off} as well as from {\em off} to {\em on}.
Recall that \eqref{eq:bestcost}, which defines a cost functional,
evaluates the jump cost at the current value of the solution pair.
For the particular case of controlling of the temperature, the jump cost should 
only depend on the current value of $q$.
Furthermore, since $q$ can only change its value at the switches, 
it needs to be forced to remain constant in between switches.
To facilitate the formulation of the optimization problem, we treat $q$ as 
an additional logic state and incorporate an input, denoted  $u \in\{0,1\}$, 
playing the role of the decision variable 
for the optimization problem.
The resulting system is given as in \eqref{eq:Hp}, 
with state $x_P = (z,q) \in \reals \times \{0,1\}$,
input $u \in \{0,1\}$, 
and data $(C_P, F_P, D_P, G_P)$ given by
\[
	\begin{aligned}
		C_P					&=\{(x_P,u):q\in\{0,1\}, u= 0\},\\
		F_P(x_P,u)	&=	(-z + z_o + z_{\Delta} q, 0)	\quad \forall (x_P,u)\in \reals^2,\\
		D_P					&=\{(x_P,u):q\in\{0,1\}, u= 1\},\\
		G_P(x_P,u) 	&=	(z,1-q)											\quad \forall (x_P,u)\in \reals^2.
	\end{aligned}
\]
With this data, flows of the plant are allowed when $u$ is zero.
In this regime, the temperature $z$ evolves according to its continuous-time
model and $q$ remains constant due $F_P$ leading to $\dot q = 0$.
At jumps, which are triggered when $u$ is equal to one, 
the update law $1 - q$ toggles the value of $q$ from $0$ to $1$ or from $1$ to $0$.

With this hybrid model, we specify 
the stage cost for flows $L_{C_P}$,
the stage cost for jumps $L_{D_P}$, 
the terminal cost $V$, and the terminal constraint set $X_P$ 
associated with Problem~\ref{prob:opti}.
As outlined above, the flow cost can be defined as an indicator of the set $\A_P := [z_{\min},z_{\max}]$ that is smooth enough.  One suitable choice is a globally Lipschitz function $L_{C_P}$ that 
depends on $z$ only and that 
in a neighborhood of $\A_P$ is equal to
\begin{equation}
L_{C_P}(z) = |z|^2_{\A},
\label{eq:lcp}
\end{equation}
while at other points has linear growth.
The jump cost is defined as a continuous function that penalizes switches.
Exploiting the fact that $q$ is a state variable, 
a suitable choice of $L_{D_P}$ that captures the cost of either transition is 
\begin{equation}
L_{D_P}(q) = 
c_{1 \to 0} q +  c_{0 \to 1} (1-q)
 \qquad \forall q \in \{0,1\},
\label{eq:ldp}
\end{equation}
where $c_{1 \to 0}$ and $c_{0 \to 1}$
are  nonnegative constants that
quantify the cost of switching the heater from 
{\em on} to {\em off} and from {\em off} to {\em on}, respectively.
The terminal cost $V$ is chosen to be equal to $L_{C_P}$, so as to quantify the distance to the desired temperature range, and the terminal constraint set $X_P$
could be simply chosen to be equal to the closed set $\A_P \times \{0,1\}$.

Next, we formulate a hybrid optimal control problem in the Mayer form in \eqref{eq:optimum} by defining the associated hybrid system $\HS$ with state $x:=(x_P,\ell)=(z,q,\ell)$, where $\ell$ is the running cost,
and its data $(C,F,D,G)$ is defined as
\[
	\begin{aligned}
		C			&:=\{x: q\in\{0,1\}\},\\ 
		F(x)	&:= (-z+z_o+z_{\Delta}q, 0, L_{C_P}(z)) \quad \forall x \in C,\\
		D			&:=C, \\ 
		G(x)	&:= (z,1-q,\ell+L_{D_P}(q))\qquad \forall x \in D.
	\end{aligned}
\]
Note that in this closed-loop formulation, since~$C=D$, jumps can occur at any time.

Given an initial condition $\xi$ for $(z,q)$ and hybrid time $(T,J) \in \realsgeq\times \nats$, 
the constraint set~$\Omega$ is chosen as
$$\Omega = \{x_0\}\times\{(T,J)\}\times X$$ 
with~$x_0:=(\xi,0)$ and~$X := X_P\times\reals = \A_P \times \{0,1\} \times \reals$,
and the cost function~$\J((z,q,\ell),(T,J),(z_\eta,q_\eta,\ell_\eta))= \ell+V(z_\eta)$.

\subsubsection{Regularity of the Optimal Control Problem}

The first question to answer is whether Problem~\ref{prob:opti} using the choices above
has a solution.  To answer this question, we apply Theorem~\ref{thrm:exist} to the Mayer 
formulation of this problem.  

Note that the associated hybrid system $\HS$ constructed above is nominally outer well-posed according to Theorem~\ref{thrm:nomout}.  
Indeed, the sets $C$ and $D$ are closed, which implies that \ref{item:A1} therein holds.
Moreover, since $L_{C_P}$ and $L_{D_P}$ are continuous by definition, the flow map $F$
and the jump map $G$ are continuous single-valued maps.  Hence, items \ref{item:A2} and
\ref{item:A3} hold.
Furthermore, for the given initial condition $\xi$ and the desired compact range of temperatures, 
the set $\Omega$ is compact.  
In addition, the cost function $\J$ is globally Lipschitz since $V = L_{C_P}$.
Finally, the hybrid system $\HS$ is such that every maximal solution is complete.  
In fact, due to the form of $C$ and $D$ combined with the regularity properties of $F$ and $G$,
\cite[Proposition 6.10]{hybridbook} implies that there exists a nontrivial solution from each initial condition 
in $C \cup D$ and that every maximal solution is complete.

Using the properties established above, by Theorem~\ref{thrm:exist}, there exists an optimal solution to the optimal
control problem \eqref{eq:optimum} if it is feasible for the given hybrid time $(T,J)$ defining $\T := \{(T,J)\}$, which implies that
Problem~\ref{prob:opti} has a solution. Moreover, by Corollary~\ref{cor:yolo0}, the optimal cost varies upper semicontinuosly.
Clearly, if $J = 0$, then, due to completeness of maximal solutions to $\HS$, there is always 
a solution. Since $C$ and $D$ do not impose any constraint on $z$, it turns
out that any choice of $(T,J) \in \realsgeq \times \nats$ leads to feasibility. These findings are summarized as follows.

\begin{theorem}
Given the cost functions~$L_{C_P}$ and~$L_{D_P}$ defined in~\eqref{eq:lcp}-\eqref{eq:ldp}, suppose that the terminal cost function~$V=L_{C_P}$ and the terminal constraint set~$X_P=[z_{\min},z_{\max}]\times\{0,1\}$. Then, Problem \ref{prob:opti} can be solved, and the optimal cost function~$h$ is upper semicontinuous.
\end{theorem}

\subsection{Bouncing Ball}

Consider a ball bouncing vertically on a horizontal flat surface, whose motion is modelled by the controlled hybrid system~$\HS_P=(C_P,F_P,D_P,G_P)$, where\footnote{The constraint~$u\in[u_{\min},u_{\max}]$ in the flow set~$C_P$ is not necessary since the flow map~$F_P$ does not depend on~$u$, but it will allow us to conclude stronger properties about inputs corresponding to optimal solution pairs.}
\[
	\begin{aligned}
		C_P					&=\{(x_P,u): p\geq 0, u\in[u_{\min},u_{\max}]\}\\
		F_P(x_P,u)	&= (v,-\gamma) \quad \forall (x_P,u)\in \reals^2\\
		D_P					&=\{ (x_P,u): p= 0, v\leq 0, u\in[u_{\min},u_{\max}]\}\\
		G_P(x_P,u)	&= (0,-\lambda v+u)	\quad \forall (x_P,u)\in \reals^2
	\end{aligned}
\]
for some~$u_{\max}\geq u_{\min}\geq 0$,~$x_P=(p,v)$ is the state with~$p\geq 0$ representing the position (height), $v$ the velocity of the ball,~$\gamma>0$ is the gravitational acceleration, and~$\lambda\in[0,1)$ is the coefficient of restitution. This system augments the canonical (autonomous) bouncing ball model with an input~$u$ that affects the post-jump velocity.

As with the thermostat, Problem~\ref{prob:opti} is equivalent with~\eqref{eq:optimum}. Given a particular instance of Problem~\ref{prob:opti} for the bouncing ball, the autonomous hybrid system~$\HS=(C,F,D,G)$ arising in this conversion to Mayer form has state~$x=(p,v,\ell)$ (with~$\ell\in\reals$ representing the running cost), where~$C=\{x: p\geq 0\}$, $D	=\{x: p= 0, v\leq 0\}$, and for every~$x\in\reals$,
\begin{equation}
	\begin{aligned}
		F(x) 	&=(v,-\gamma,0,L_{C'}(p,v)),\\
		G(x) 	&=\{(0,-\lambda v+u,\ell+L_{D'}(v,u)): u\in[u_{\min},u_{\max}]\}.
	\end{aligned}
\label{eq:equation}
\end{equation}

The cost function~$L_{C'}$ above is introduced to simplify the problem and replace~$L_{C_P}$, as the flow map~$F_P$ does not depend on~$u$. Similarly, the cost function $L_{D'}$ is introduced as~$p=0$ at jumps. The constraint set of the problem is given as~$\Omega = \{x_0\}\times\{(T,J)\}\times X$ with~$x_0=(\xi,0)$ and~$X= X_P\times\reals$, and the cost function~$\J(x,(T,J),\eta)= \ell+V(p,v)$, where~$X_P$ and~$V$ are the terminal constraint set and terminal cost function in Problem~\ref{prob:opti}.

\subsubsection{Well-Posedness, Existence, and Upper Semicontinuity}

We show that when the cost functions~$L_{C'}$ and~$L_{D'}$ are continuous (on~$C':=\{(p,v): p\geq 0\}$ and~$D':=\{(v,u): v\leq 0, u\in[u_{\min},u_{\max}]\}$, respectively), the augmented system~$\HS$ is nominally inner and outer well-posed. Nominal outer well-posedness follows directly from Theorem~\ref{thrm:nomout}. Indeed, the sets~$C$ and~$D$ are closed, \ref{item:A2} is due to the flow map~$F$ being single-valued and continuous on the flow set~$C$, and~\ref{item:A3} is due to the jump map~$G$ being compact-valued and continuous on the jump set~$D$.

For nominal inner well-posedness, let
\[
\widehat{G}(x):=(0,-\lambda v,\ell+L_{D'}(v,u_{\min})) \quad\forall x\in\reals^3.
\]
Consider the hybrid system~$(C,F,D,\widehat{G})$ and note that maximal solutions of this system are complete and unique, due to the following reasons: a) the first two components of a given solution corresponds to a solution of the autonomous system arising from~$\HS_P$ when the input~$u=u_{\min}$, which is bounded, b) when~$u=u_{\min}$ for $\HS_P$, the resulting autonomous system has unique maximal solutions, and c)~$L_{C'}$ is continuous on~$C':=\{(p,v):h\geq 0\}$, which implies integrability. Then, by \cite[Proposition 7]{arxiv},~$(C,F,D,\widehat{G})$ is nominally inner well-posed, which implies the continuous-time system~$\dot{x}=F(x)$ $x\in C$ is nominally inner well-posed. From \cite[Theorem 17 and Proposition 19]{arxiv}, it follows that the hybrid system~$\HS$ is nominally inner well-posed if the jump map $G$ and the mapping~$\widetilde{G}$ given below are inner semicontinuous relative to the jump set~$D$:
\[
	\widetilde{G}(x):= G(x)\cap (\widetilde{C}\cup D) \quad \forall x\in\reals^3,
\]
where~$\widetilde{C}\subset \closure C$ is the set of points from where the constrained differential equation~$\dot{x}= F(x)$ $x\in C$ has nontrivial solutions. Hence, it follows that~$\widetilde{C}\cup D= C$. Since~$G$ is inner semicontinuous and~$G(x)\cap C \subset C$ for all~$x\in D$, it follows that~$\HS$ is nominally inner well-posed. Finally, we note that~$\HS$ is pre-forward complete, which is due to the fact that the constrained differential equation~$\dot{x}= F(x)$ $x\in C$ has bounded solutions.

Thus, due to the equivalence of Problem~\ref{prob:opti} and~\eqref{eq:optimum}, the following can be concluded.
\begin{thm}
Suppose that the cost functions~$L_{C'}$ and~$L_{D'}$ are continuous on the sets~$C'=\{(p,v):h\geq 0\}$ and~$D'=\{(v,u): v\leq 0, u\in[u_{\min},u_{\max}]\}$, respectively. Then, if the constraint set~$X_P$ is closed and the cost function~$V$ is lower semicontinuous on~$X_P$, the optimal cost function~$h$ is upper semicontinuous. If, in addition, there exists a solution pair~$(x_P,u)$ such that~$x_P(0,0)=\xi$, $(T,J)\in \dom (x_P,u)$, and~$x_P(T,J)\in X_P$, then Problem \ref{prob:opti} can be solved.
\end{thm}
\begin{pf}
This can be directly inferred via Theorems~\ref{thrm:exist} and \ref{thrm:usc} by replacing the set~$X= X_P\times\reals$ with its intersection with the reachable set~$\reach_{\HS}(x_0,T,J)$, which is compact by \cite[Proposition~4.2]{cdc2020}.
\end{pf}

\subsubsection{Continuity of the Optimal Cost and Outer/Upper Semicontinuity of Optimal Solutions}

To conclude stronger properties about Problem~\ref{prob:opti}, we assume that the terminal cost~$V$ is continuous on the set~$C'$, which would imply that the resulting closed-loop cost function~$V'$ is continuous on~$\closure(C) \cup D$, and the terminal constraint set~$X_P=\reals^n$. In the sequel, we rely on Corollaries~\ref{cor:gme} and~\ref{cor:coupdegrace}. To invoke these corollaries, it is necessary and sufficient to also show that the reachable set~$\reach_{\HS}(x_0,T,J)$ is nonempty, and for every~$\xi\in\reach_{\HS}(x_0,T,J)$, there exists~$x\in\sol_{\HS}(x_0)$ such that~$\xi=x(T,J)$ and~$T$ is not a jump time or the terminal ordinary time of~$x$.

Given the initial condition~$\xi = (\xi_1,\xi_2)$ and a parameter~$\nu\in[u_{\min}, u_{\max}]$, let~$(x,u)$ be the unique solution pair with~$x(0,0)\in \xi$ and~$\dom(x,u)$ unbounded, satisfying~$u(t,j)=u(s,i)=\nu$ for all~$(t,j),(s,i)\in\dom (x,u)$. Existence of such a pair is easy to show following an analysis similar to the one in the previous subsection. Regardless of the choice of the input~$\nu$, the first impact with the ground occurs at ordinary time~$(\xi_2+\sqrt{\xi_2^2+2\gamma\xi_1})/\gamma$ (\cite[Example~2.12]{hybridbook}) with velocity~$-\sqrt{\xi_2^2+2\gamma\xi_1}$---the latter can be derived using conservation of energy during flows. Given~$j\geq 1$, let~$t_j$ be the ordinary time of jump~$j$ and let~$v_j\geq 0$ be the velocity of the ball immediately after jump~$j$, i.e., $v_j:=x(t_j,j)$. Then,~$v_1=\lambda\sqrt{\xi_2^2+2\gamma\xi_1}+\nu$. Moreover,~$t_{j+1}-t_j=2v_j/\gamma$ by \cite[Example~2.12]{hybridbook} and~$v_{j+1}=\lambda v_j+\nu \geq 0$ (again due to conservation of energy, which implies that the velocity right after jump~$j$ and right before~$j+1$ differ only in their sign) for all~$j\geq 1$. From these equations, one can then derive
\begin{equation}
\label{eq:jumptimes}
	\begin{aligned}
		t_1 				&= \left(\xi_2+\sqrt{\xi_2^2+2\gamma\xi_1}\right)/\gamma,\\
		t_{j}^{\nu}	&= t_1+\frac{2}{\gamma(1-\lambda)}\times \psi(\nu,j) \quad \forall j\geq 2,
	\end{aligned}
\end{equation}
where the superscript~$\nu$ is included to indicate dependency on the input parameter~$\nu$, and
\begin{multline}
\label{eq:psi}
\psi(\nu,j):=(j-1)\nu +\\
\left(\underbrace{\lambda\sqrt{\xi_2^2+2\gamma\xi_1}+\nu}_{v_1} -\frac{\nu}{1-\lambda}\right)(1-\lambda^{j-1})
\end{multline}
with~$v_1$ indicating (post-impact) velocity after the first jump, i.e.~$x(t_1,1)$. Since~$t_j^{\nu}u$ is an increasing function of~$\nu$ for fixed~$j$, one can then infer the following: a) the reachable set~$\reach_{\HS}(x_0,T,J)$ of the augmented autonomous system is nonempty if~$t_{J}^{u_{\min}}\leq T\leq t_{J+1}^{u_{\max}}$, b) for every~$\xi\in\reach_{\HS}(x_0,T,J)$ there exists~$x\in\sol_{\HS}(x_0)$ such that~$\xi=x(T,J)$ and~$T$ is not a jump time or the terminal ordinary time of~$x$ if~$t_{J}^{u_{\max}}\leq T\leq t_{J+1}^{u_{\min}}$.

Using these two deductions, we reach the following result, where $\mathcal{O}_{\HS}(\xi,T,J)$ denotes the set of optimal solution pairs of Problem~\ref{prob:opti}, that is the set of solution pairs~$(x_P,u)$ of~$\HS_P$ subject to the constraints~$x_P(0,0)=\xi$ and~$x_P(T,J)\in X_P$.

\begin{thm}
\label{thrm:yesbest}
Suppose that the cost functions~$L_{C'}$ and~$L_{D'}$ are continuous on the sets~$C'=\{(p,v):h\geq 0\}$ and~$D'=\{(v,u): v\leq 0, u\in[u_{\min},u_{\max}]\}$, respectively. Moroever, suppose that the terminal constraint set~$X_P=\reals^n$, the terminal cost~$V$ is continuous on~$C'$, and~$t_{J}^{u_{\max}}\leq T\leq t_{J+1}^{u_{\min}}$, where~$t_{j}^{\nu}$ is defined in~\eqref{eq:jumptimes}-\eqref{eq:psi}. Then, Problem \ref{prob:opti} can be solved, and the function~$h$ is continuous at~$(\xi,T,J)$. Moreover, the following hold.
\begin{description}
	\item[Local Boundedness] There exists~$\varepsilon>0$ and a compact set~$K$ such that
	\begin{multline*}
		\xi'\in\xi+\varepsilon\ball,\,(T',J')\in(T,J)+\varepsilon\ball\\
		\implies (x_P(t,j),u(t,j))\in K \quad
	\end{multline*}
	for all $(t,j)\in\dom (x,u)$ and $(x_P,u)\in\mathcal{O}_{\HS}(\xi',T',J')$.
	\item[Outer Semicontinuity] Let~$\{(x'_{P_i},u'_i)\}_{i=0}^{\infty}$ be a sequence of optimal solution pairs such that $(x'_{P_i},u'_i)\in\mathcal{O}_{\HS}(\xi_i,T_i,J_i)$ for all~$i\geq 0$ and~$\{x'_{P_i}\}_{i=0}^{\infty}$ is graphically convergent. Then, if the sequences~$\{\xi_i\}_{i=0}^{\infty}$ and~$\{(T_i,J_i)\}_{i=0}^{\infty}$ converge to~$\xi$ and~$(T,J)$, respectively, there exists~$u$ such that~$(x_P,u)\in\mathcal{O}_{\HS}(\xi,T,J)$,	where~$x_P$ is the graphical limit of~$\{x'_{P_i}\}_{i=0}^{\infty}$. In addition, given~$j\geq 0$, if~$t_{j+1}$ is the $(j+1)$-th jump time of~$(x_P,u)$, then~$u((t_{j+1},j))=\lim_{i\to\infty}u'_i(t^i_{j+1},j)$, where~$t^i_{j+1}$ is the~$(j+1)$-th jump time of~$(x'_{P_i},u'_i)$ for large enough~$i$.
	\item[Upper Semicontinuity] For all~$\tau\geq 0$ and~$\varepsilon>0$, there exists~$\eta>0$ such that the following holds: for every $x'_0\in x_0+\eta\ball$, $(T',J')\in(T,J)+\eta\ball$, and~$(x'_P,u')\in\mathcal{O}_{\HS}(x'_0,T',J')$, there exists~$(x_P,u)\in\mathcal{O}_{\HS}(x_0,T,J)$ such that~$x$ and~$x'$ are $(\tau,\varepsilon)$-close, and the following holds for every~$j\geq 0$: if~$t_{j+1}$ and~$t'_{j+1}$ are the $(j+1)$-th jump times of~$(x,u)$ and~$(x',u')$, respectively, then~$|u(t_{j+1},j)-u'(t'_{j+1},j)|<\varepsilon$.
\end{description}
\end{thm}
\begin{pf}
Existence of an optimal solution pair and continuity of~$h$ follows directly from Corollary~\ref{cor:gme}. Local boundedness is due to the local boundedness result in Corollary~\ref{cor:coupdegrace} and the fact that the inputs are constrained to the compact set~$[u_{\min},u_{\max}]$. For the second item, graphical convergence of the state trajectories to~$x$ is proved in Corollary~\ref{cor:coupdegrace}, and the statement regarding the inputs at jump times follows from graphical convergence of the state trajectories, \cite[Lemma~2]{arxiv}, continuity of~$G_P$, and the fact that the mapping~$G_P(x_P,.)$ is one-to-one. The final statement can then be proven by contradiction as in the proof of Theorem~\ref{thrm:coupdegrace}.
\end{pf}

\subsubsection{Illustrative Example}

To numerically illustrate the results, we consider the control problem in~\cite{9147972} of ensuring that the ball reaches a desired peak height~$p_{{\rm des}}$ after every impact, asymptotically as the number of jumps tends to infinity. Equivalently, due to conservation of energy during flows, the control objective can be viewed as asymptotically stabilizing the set~$\A:=\{x_P: W(x_P)=W(p_{{\rm des}})\}$, where $W(x_P):=\gamma p+ v^2/2$ is the total energy function, which is continuous. To achieve this objective, we let the terminal cost function~$V=W$, and select~$L_{C'}$ as the zero function. For the jump cost, let~$L_{D'}(v,u)=\gamma p_{{\rm des}}(v+\sqrt{2\gamma p_{{\rm des}}})^2/2$ if~$ v\geq -\sqrt{2\gamma p_{{\rm des}}}/\lambda$, otherwise, let
\begin{multline*}
	L_{D'}(v,u)=\min\Big\{\gamma p_{{\rm des}}(v+\sqrt{2\gamma p_{{\rm des}}})^2/2,\\
	(v^2/2-\gamma p_{{\rm des}})^2-(\lambda^2v^2/2-\gamma p_{{\rm des}})^2\Big\},
\end{multline*}
which is continuous.

\begin{figure}
\centering
\subfloat[]{
\label{fig:voltage}
\includegraphics[width=1.0\columnwidth]{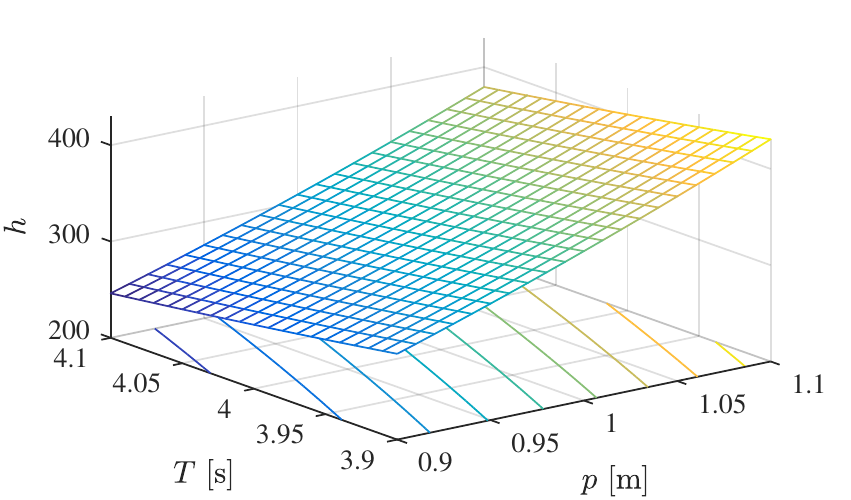}}
\vspace{8pt}
\subfloat[]{
\label{fig:phaseplot}
\includegraphics[width=1.0\columnwidth]{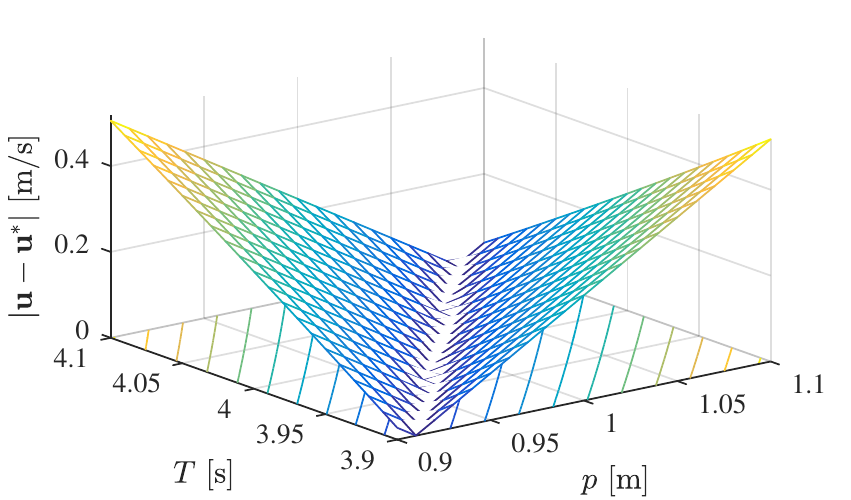}}
\caption[]{Simulation results for the bouncing ball with $J = 2$ and initial velocity~$v=0$, as the ordinary time horizon parameter~$T$ and the initial height~$p$ are varied. \subref{fig:voltage} Optimal cost as~$T$ and~$p$ are varied. \subref{fig:phaseplot} Convergence of the optimal input as~$T$ and~$p$ tends to the nominal values of~$T=4$ and~$p=1$; the vector $\mathbf{u}\in\reals^2$ represents the values of the optimal input at jump times, $\mathbf{u}^*\in\reals^2$ corresponds to the case~$T=4$ and~$p=1$.}
\label{fig:ex3}%
\end{figure}

\begin{figure}[tbhp]
\centering
	\includegraphics[width=1.0\columnwidth]{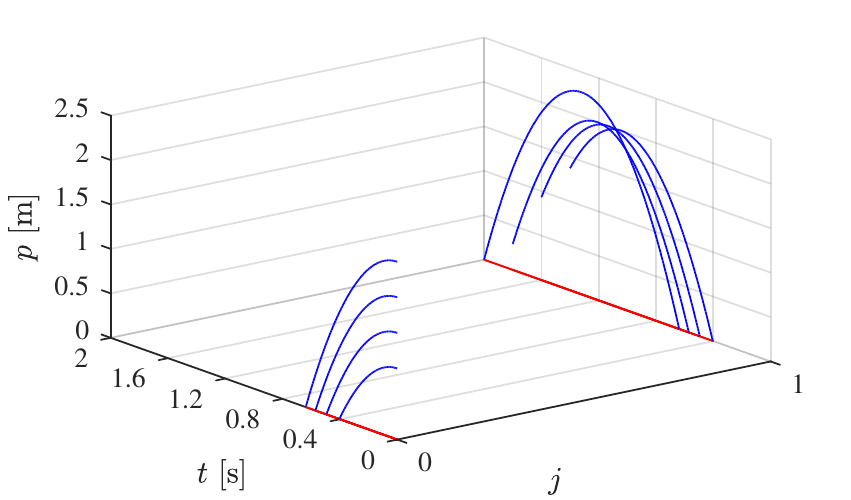}
\caption{Graphical convergence of the state trajectories with $J = 1$ and initial velocity~$v=0$, as the ordinary time horizon parameter~$T$ and the initial height~$p$ tend to nominal values of~$T=2$ and~$p=2$.}
\label{fig:observer}
\end{figure}

Simulation results\footnote{Code can be found at \url{https://github.com/HybridSystemsLab/HybridOptimalControlBouncingBall}} corresponding the parameters $\gamma = 9.81$\,m/s$^2$, $\lambda = 0.8$, $u_{\min} = 1$\,m/s, $u_{\max} = 10$\,m/s, and $p_{\rm des}=\,2$ m can be seen in Figures \ref{fig:ex3} and \ref{fig:observer}. The optimal control problem is solved by casting it as a nonlinear program with linear inequalities, using the closed-form analytical solutions of the system and the fact that the flow cost function~$L_{C'}$ is zero. Although the condition regarding jump times in Theorem \ref{thrm:yesbest} are not verified explicitly\footnote{This can lead to conservative results when~$u_{\max}-u_{\min}$ is large.}, the findings of the theorem regarding continuity of the optimal cost, graphical convergence of the trajectories, and convergence of the inputs at jump times are still observed. In particular, Figure \ref{fig:ex3} shows that the optimal cost depends continuously on the initial height and the ordinary time horizon in a neighborhood of their nominal values, and similarly, the optimal input (at jump times) depends continuously on the initial height and the ordinary time horizon at their nominal values. In Figure \ref{fig:observer}, outer semicontinuous dependence of the optimal state trajectories on the initial height and the ordinary time horizon can be observed: as the initial height and the time horizon parameter converge to their nominal values, the corresponding optimal state trajectories similarly converge (graphically) to the optimal state trajectory corresponding to these nominal values.

\section{Conclusion}
\label{sec:conc}
Existence of optimal solutions, their dependency on constrains
and perturbations, as well as properties revealing the dependency of the optimal cost
and of the value function with respect to the given data
and perturbations are established for a general hybrid control problem.
It is shown that nominal outer well-posedness of the hybrid dynamical system
is instrumental guaranteeing not only the existence of an optimal 
solution to the hybrid optimal control problem but also upper semicontinuity of the value function, plus its upper semicontinuous dependence  on initial conditions and perturbations.  
In addition, when the hybrid dynamical system is inner well-posed,
the optimal cost is continuous and, very importantly, can be continuously approximated.  

With sufficient conditions for outer and inner well-posedness 
for the class of hybrid systems considered already being available in the literature, the results in this paper pave the road to the understanding of the effect of computation and approximation in emerging tools for hybrid dynamical systems, such as numerical simulation, model predictive control, and parameter estimation. Future work includes exploiting the continuous approximation of the optimal cost to the model predictive control framework proposed in \cite{ALTIN2018128,acc,9147972} when the hybrid system to control is discretized for the purpose of solving the optimal control problem associated with model predictive control.


%


\bibliographystyle{unsrt}
\bibliography{bibs/reach,bibs/long,bibs/Biblio,bibs/RGSweb}


\appendix

\section{Closeness of Hybrid Arcs}
\label{sec:extras}

The following recalls~\cite[Definition~5.23]{hybridbook}.
\begin{defn}
Given~$\tau\geq 0$ and $\varepsilon>0$, two hybrid arcs~$x$ and~$x'$ are said to be~$(\tau,\varepsilon)$-close if
	\begin{itemize}[leftmargin=*]
		\item for every~$(t,j)\in\dom x$ satisfying~$t+j\leq\tau$, there exists~$(t',j)\in\dom x'$ such that $|t-t'|<\varepsilon$ and $|x(t,j)-x'(t',j)|<\varepsilon$;
		\item for every~$(t',j')\in\dom x'$ satisfying~$t'+j'\leq\tau$, there exists~$(t,j')\in\dom x$ such that $|t'-t|<\varepsilon$ and $|x'(t',j')-x(t,j')|<\varepsilon$.
	\end{itemize}
\end{defn}

\section{Sufficient Conditions for Nominal Well-Posedness}
\label{sec:hbc}

Although (nominal) outer well-posedness of a hybrid system can be difficult to check, it is guaranteed when the data of the system satisfies the so-called hybrid basic conditions~\cite[Theorem~6.8]{hybridbook}.

\begin{thm}
\label{thrm:nomout}
A hybrid system~$\HS=(C,F,D,G)$ is (nominally) outer well-posed if the following hold.
\begin{enumerate}[label={(A\arabic*)},leftmargin=*]
	\item \label{item:A1}	The sets~$C$ and~$D$ are closed.
	\item \label{item:A2}	The flow map~$F$ is locally bounded and outer semicontinuous relative to~$C$, and~$C\subset \dom F$. Furthermore, for every~${x\in C}$, the set~${F(x)}$ is convex.
	\item \label{item:A3}	The jump map~$G$ is locally bounded and outer semicontinuous relative to~$D$, and~$D\subset \dom G$.
\end{enumerate}
\end{thm}

Given a set~$S\subset\reals^n$ and a point~$x\in\reals^n$, denote by~$T_S(x)$ the \textit{Bouligand tangent cone} to~$S$ at~$x$ \cite[Definition 5.12]{hybridbook} and by~$M_S(x)$ the \textit{Dubovitsky-Miliutin tangent cone} to~$S$ at~$x$ \cite[Definition 4.3.1]{viability}. A set of sufficient conditions for nominal inner well-posedness, which use these tangent cones, are given below~\cite[Theorem~1.1]{cdc2020}. For a proof of this result, see the discussion at the end of \cite[Section 5.2]{arxiv}.

\begin{thm}
\label{thrm:viab}
Given a hybrid system~$\HS = (C,F,D,G)$, suppose that the flow set~$C$ is closed and~\ref{item:A2} holds. Then,~$\HS$ is nominally inner well-posed if the following hold.
\begin{enumerate}[label={(B\arabic*)},leftmargin=*]
	\item \label{item:R1}	For every~$x\in C$, there exists an extension of~$F|_{C}$ that is closed valued and Lipschitz\footnote{A set-valued mapping $M$ is \textit{Lipschitz} on~$X$ if it has nonempty values on~$X$ and there exists~$L\geq 0$ such that~$M(x)\subset M(x')+L|x-x'|\ball$ for every~$x,x'\in X$.} on a neighborhood of~$x$.
	\item \label{item:R2}	For every~$x\in \partial C$ such that~$F(x)\cap T_C(x)$ is nonempty, there exists~$r>0$ such that~$F(x')\subset M_{\interior C}(x')$ for all~$x'\in(x+r\ball)\cap\partial C$, and~$(x+r\ball)\cap D \subset C$.
	\item \label{item:R3} For every~$x\in\interior C\cap\partial D$,~$F(x)\cap M_{\interior D}(x)$ is nonempty.
	\item \label{item:R4} For every~$x\in\partial C\cap \partial D$, either of the following hold:
	\begin{itemize}
		\item there exists~$r>0$ such that $(x+r\ball)\cap C\subset D$;
		\item $F(x)\cap M_{\interior C}(x)\cap M_{ \interior D}(x)$ is nonempty;
		\item $F(x)\cap T_C(x)$ is empty and there exists~$r>0$ such that $(x+r\ball)\cap \partial C\subset D$.
	\end{itemize}
	\item \label{item:N5} The jump map~$G$ is inner semicontinuous relative to~$D$.
	\item \label{item:N6} The mapping~$\widetilde{G}:\reals^n\rightrightarrows\reals^n$, where
	\begin{equation*}
		\begin{aligned}
			\widetilde{G}(x)&:=G(x)\cap(\widetilde{C}\cup D) \quad \forall x\in \reals^n,\\
			\widetilde{C}							&:=\interior(C)\cup\{x\in\partial C: F(x)\cap T_C(x)\neq\varnothing \},
		\end{aligned}
	\end{equation*}
	is inner semicontinuous relative to~$D$.
\end{enumerate}
\end{thm}

\end{document}